\documentclass{svmult}

\usepackage{makeidx}         
\usepackage{graphicx}        
\usepackage{multicol}        

\makeindex             
\begin{document}

\title*{Quadratic optimal functional quantization of
stochastic processes and  numerical applications}
\titlerunning{Functional Quantization}
\author{Gilles Pag\`es\inst{1}}
\institute{Laboratoire de Probabilit\'es et Mod\`eles al\'eatoires, UMR~7599, Universit\'e Paris 6, case 188, 4,
pl. Jussieu, F-75252 Paris Cedex 5, France.
\texttt{gpa@ccr.jussieu.fr}
}
%
%

\maketitle

\begin{abstract}
In this paper, we present an overview of the recent developments of functional 
quantization of stochastic processes, with an emphasis on the quadratic case. Functional quantization  is a way to
approximate a process, viewed as a Hilbert-valued random variable, using a nearest neighbour projection on a finite
codebook.   A special emphasis is made on the computational aspects and the numerical applications, in particular the pricing
of  some path-dependent European options.
\end{abstract}

\printindex



\font\tenmath=msbm9 scaled 1200
\font\sevenmath=msbm7 scaled 1200
\font\fivemath=msbm5 scaled 1200

\newfam\mathfam \textfont\mathfam=\tenmath
\scriptfont\mathfam=\sevenmath \scriptscriptfont\mathfam=\fivemath
\def\math{\fam\mathfam}
\def\R{{\math R}}
\def\N{{\math N}}
\def\E{{\math E}}
\def\L{{\cal L}}
\def\P{{\math P}}
\def\Z{{\math Z}}
\def\T{{\math T}}
\def\Q{{\math Q}}
\def\B{{\math B}}
\def\C{{\math C}}
\def\sko{{\math D}}
\def\F{{\cal F}}

%
%

\section{Introduction}
Functional quantization  is a way to discretize the path space of a
stochastic process. It has been extensively investigated since the
early 2000's by several authors (see among
others~\cite{LUPA1},~\cite{LUPA3},~\cite{Dereichetal},
\cite{DEFEMASC},~\cite{LUPA2}, etc). It first appeared as a natural
extension of the   Optimal Vector Quantization theory of
(finite-dimensional)  random vectors which finds its origin in the
early 1950's for signal processing (see~\cite{GEGR}
or~\cite{GRLU1}).

Let us consider   a  Hilbertian setting. One considers
a  random vector $X$  defined on a probability space $(\Omega,{\cal A},\P)$ taking its values
in a separable Hilbert space $(H, (.|.)_{_H})$ (equipped with its natural Borel $\sigma$-algebra) and satisfying
$\E|X|^2<+\infty$.  When $H$ is an Euclidean space ($\R^d$), one speaks about {\em
Vector Quantization}. When $H$ is an infinite dimensional space like
$L^2_{_T}:= L^2([0,T], dt)$ (endowed with the usual Hilbertian norm
$|f|_{L^2_{_T}}:= (\int_0^Tf^2(t)dt )^{\frac12}$)  one
speaks of {\em functional quantization} (denoted  $L^2_{_T}$ from
now on). A (bi-measurable) stochastic process $(X_t)_{t\in [0,T]}$
defined on $(\Omega,{\cal A},\P)$ satisfying
$|X(\omega)|_{L^2_{_T}}<+\infty$ $\P(d\omega)$-$a.s.$  can always be
seen, once possibly modified on a $\P$-negligible set,  as an
$L^2_{_T}$-valued random variable. Although we will focus on the
Hilbertian framework, other choices are possible for $H$, in
particular some more general Banach settings like $L^p([0,T],dt)$ or
${\cal C}([0,T],\R)$ spaces.

This paper is organized as follows: in Sections~\ref{Defquadraquant} we introduce  quadratic quantization in a
Hilbertian setting. In Section~\ref{Optquadratic}, we focus on optimal quantization, including some extensions to non
quadratic quantization.  Section~\ref{Cubature} is devoted to some quantized cubature formulae. Section~\ref{RateRd} provides
some classical background on the quantization rate in finite dimension. Section~\ref{OptiQuantB}  deals with
functional quantizations of Gaussian processes, like the Brownian motion, with a special emphasis on the numerical
aspects. We present here what is, to our guess,  the first large scale  numerical optimization of the quadratic quantization of
the Brownian motion. We compare it to the optimal product quantization, formerly investigated in~\cite{PAPR2}. In section, we
propose a constructive approach to the functional quantization of scalar or multidimensional  diffusions (in the
Stratanovich sense). In Section~\ref{pathdeppric}, we show how to use  functional quantization to price path-dependent
options like Asian options (in a heston stochastic volatility model). We conclude  by   some recent results showing
how to derive universal (often optimal) functional quantization rate from time regularity of a process
in Section~\ref{FQReg} and by a few clues in Section~\ref{LowBounds}  about the specific methods that
produce some lower bounds (this important subject as many others like the connections with small deviation theory  is not
treated in this numerically oriented overview. As concerns statistical applications of functional quantization we refer
to~\cite{TAKI, TAPEOG}.

\medskip
\noindent  {\sc Notations.} $\bullet$ $a_n\approx b_n$ means $a_n=O(b_n)$ and $b_n=O(a_n)$; $a_n\sim b_n$ means $a_n=
b_n +o(a_n)$.

\smallskip
\noindent $\bullet$ If $X:(\Omega, {\cal A}, \P)\to (H,|\,.\,|_{_H})$ (Hilbert space), then $\|X\|_{_2}=
(\E|X|_{_H}^2)^{\frac 12}$.

\smallskip
\noindent $\bullet$ $\lfloor x\rfloor$ denotes the integral part of the real $x$.

\section{What is quadratic functional quantization?}\label{Defquadraquant}
\begin{figure}\label{fig:1}
\centering
\begin{tabular}{c}
 \centerline{\includegraphics[width=7cm,height = 4.5 cm]{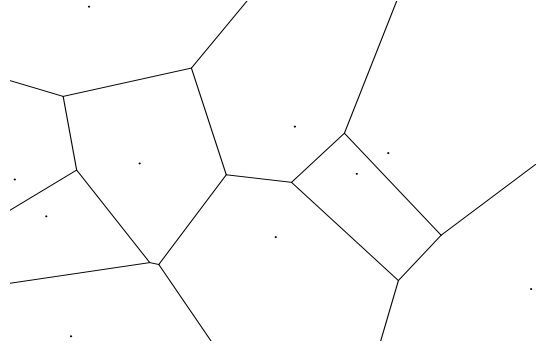}}
\end{tabular}
\caption{\em A $2$-dimensional $10$-quantizer
$\Gamma=\{x_1,\dots,x_{10}\}$ and its Voronoi diagram.}
\end{figure}

Let $(H,(\,.|.\,)_{_H})$ denote a separable Hilbert space. Let $X\!\in L^2_{_H}(\P)$ $i.e.$  a random vector
$X:(\Omega, {\cal A}, \P)\longmapsto H$ ($H$ is endowed with its Borel $\sigma$-algebra) such that
$\E\,|X|^2_{_H}<+\infty$.  An {\em $N$-quantizer} (or {\em $N$-codebook}) is  defined as a subset
\[
 \Gamma:=\{x_1,\ldots,x_{_N}\}\subset  H
\]
with card$\Gamma=N$. In  numerical applications, $\Gamma$ is also called  {\em grid}. Then, one can {\em quantize}
(or simply  discretize) $X$ by $q(X)$ where $q:H\mapsto \Gamma$ is a Borel function. It is
straightforward that
\[
\forall\, \omega\!\in \Omega,\qquad |X(\omega)-q(X(\omega))|_{_H}\ge d(X(\omega), \Gamma)= \min_{1\le i\le
N}|X(\omega)-x_i|_{_H}
\]
so that the best pointwise approximation of $X$ is provided by considering for  $q$ a nearest neighbour projection on
$\Gamma$, denoted ${\rm Proj}_{_\Gamma}$. Such a projection is in one-to-one correspondence
with the Voronoi partitions (or diagrams) of $H$ induced by $\Gamma$ $i.e.$  the  Borel partitions of $H$ satisfying
\[
C_i(\Gamma)\subset\left\{\xi\!\in H: |\xi-{x}_i|_{_H} = \min_{1\le j\le
N}|\xi-{x}_j|_{_H}\right\}=\overline{C}_i(\Gamma),\qquad i=1,\ldots,N,
\]
where $\overline{C}_i(\Gamma)$ denotes the closure of $C_i(\Gamma)$ in $H$ (this heavily uses the Hilbert structure). Then
$$
{\rm Proj}_{_\Gamma}(\xi) :=
\sum_{i=1}^N x_i
\mbox{\bf 1}_{C_i(\Gamma)}(\xi)
$$
is a nearest neighbour projection on $\Gamma$. These projections  only differ
on the boundaries of the  {\em Voronoi cells}
$C_i(\Gamma)$, $i=1,\ldots,N$. All Voronoi partitions have the same boundary contained in the    union of
the median hyperplanes defined by the pairs $(x_i,x_j)$, $i\!\neq\!j$. Figure~\ref{fig:1} represents the Voronoi
diagram defined by a (random) $10$-tuple in $\R^2$.
Then, one defines a {\em Voronoi $N$-quantization} of $X$ by setting for every
$\omega\!\in\Omega$,
\[
\widehat{X}^{\Gamma}(\omega):= {\rm Proj}_{ \Gamma}(X(\omega))= \sum_{i=1}^N x_i \mbox{\bf
1}_{C_i(\Gamma)}(X(\omega)).
\]
One clearly has, still for every $\omega\!\in\Omega$, that
\[
|X(\omega)-\widehat{X}^{\Gamma} (\omega)|_{_H}  = {\rm
dist}_{_H}(X(\omega),\Gamma) = \min_{1\le i\le N}
|X(\omega)-x_i|_{_H}.
\]
The mean (quadratic) quantization error is then defined by
\begin{equation}\label{eNGammaXH}
e(\Gamma,X,H)= \|X-\widehat{X}^{\Gamma}\|_{_2} = \sqrt{\E
\left(\min_{1\le i\le N} |X-{x}_i|_{_H} ^2 \right)}.
\end{equation}
The distribution of $\widehat X^{\Gamma}$ as a random vector is given by the $N$-tuple $(\P(X\!\in
C_i(\Gamma)))_{1\le i\le N}$ of the Voronoi cells. This distribution clearly depends on the choice of the
Voronoi partition as emphasized by the following elementary situation: if $H=\R$, the distribution of $X$
is given by $\P_{_X}= \frac 13(\delta_0+\delta_{1/2}+\delta_1)$, $N=2$ and $\Gamma=\{0,1\}$ since
$1/2\in\partial C_0(\Gamma)\!\cap\!\partial C_1(\Gamma)$. However, if
$\P_{_X}$ weights no hyperplane, the distribution of $\widehat X^{\Gamma}$ depends only on $\Gamma$.

\medskip   As concerns terminology, {\em Vector Quantization} is concerned with the finite dimensional case
--~when ${\rm dim}H<+\infty$~-- and is a rather old story, going back to the early 1950's when it was designed in
the field of signal processing and then mainly developed in the community of Information Theory.
The term  {\em functional quantization}, probably introduced in~\cite{PAG0.5, LUPA1}, deals with the
infinite dimensional case including the more   general Banach-valued setting. The term
``functional" comes from the fact  that a typical infinite dimensional Hilbert space is
the function space $H=L^2_{_T}$. Then, any   (bi-measurable) process $X:([0,T]\times
\Omega,Bor([0,T])\otimes {\cal A})\to
(\R,Bor(\R))$ can be seen as a random vector taking values in the set of Borel functions on $[0,T]$. Furthermore,
$((t,\omega)\mapsto X_t(\omega))\!\in L^2(dt\otimes  d\P )$ if and only if $(\omega \mapsto X_.(\omega))\!\in
L^2_H(\P)$ since
\[
\int_{[0,T]\times \Omega} X^2_t(\omega)\, dt\,\P(d\omega)=
\int_{\Omega} \P(d\omega)\int_0^T X^2_t(\omega) \,dt=
\E\,|X_.|^2_{L^2_{_T}}.
\]

\section{Optimal (quadratic) quantization}\label{Optquadratic}
At this stage we are lead  to wonder whether it is possible to design some optimally fitted  grids   to a given
distribution $\P_{_X}$ $i.e.$ which induce the lowest possible mean quantization error among all grids of
size at most $N$. This amounts to the following optimization problem
\begin{equation}
e_{_{N}}(X,H) := \inf_{\Gamma\subset H, \mbox{\footnotesize card}(\Gamma) \le N}e(\Gamma,X,H).
\end{equation}
\begin{center}
\begin{figure}
\begin{tabular}{ll}
\includegraphics[width=6cm,height = 5.25cm]{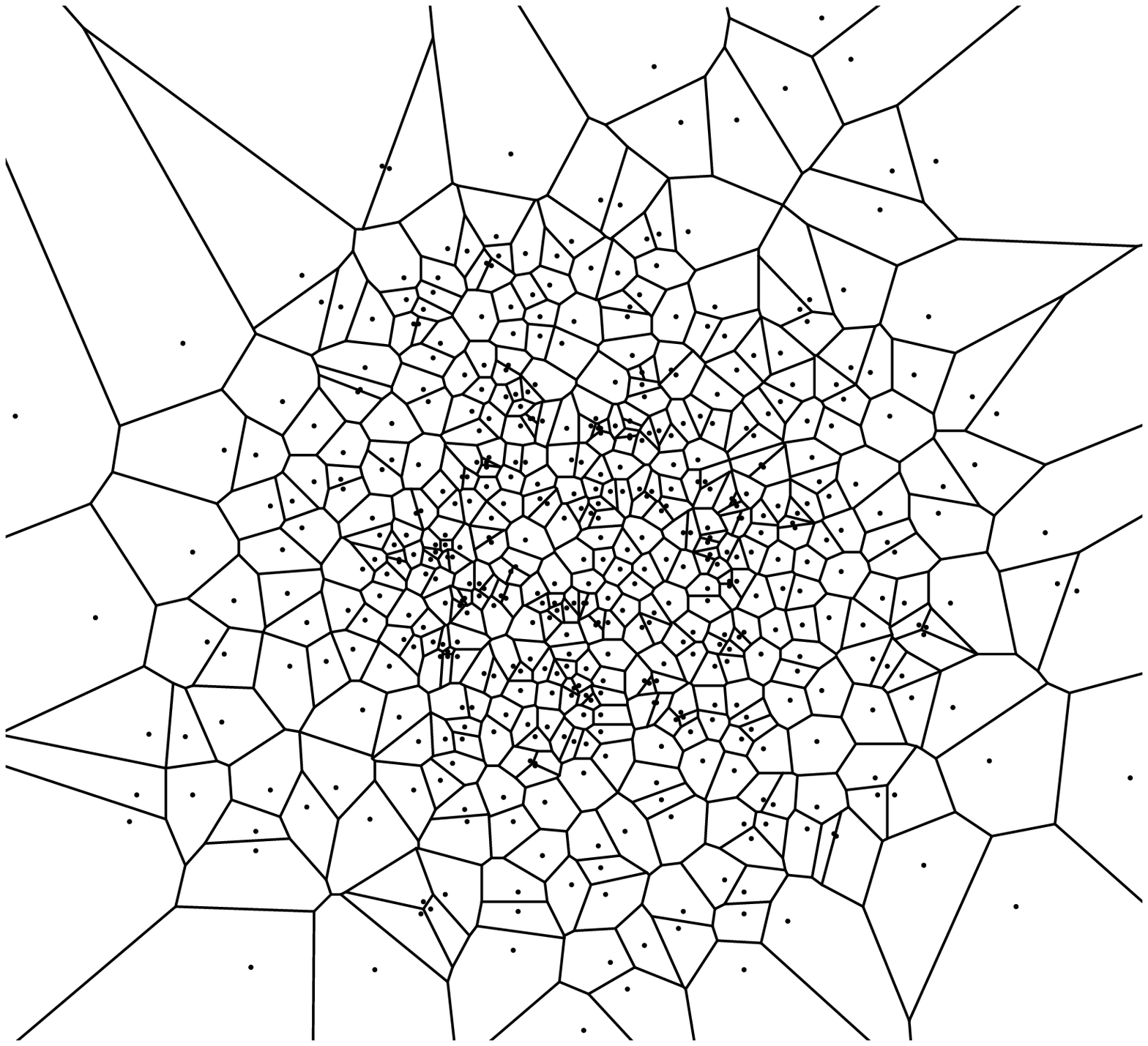}$\quad$&\includegraphics[width=5.25cm,height =
5.25cm]{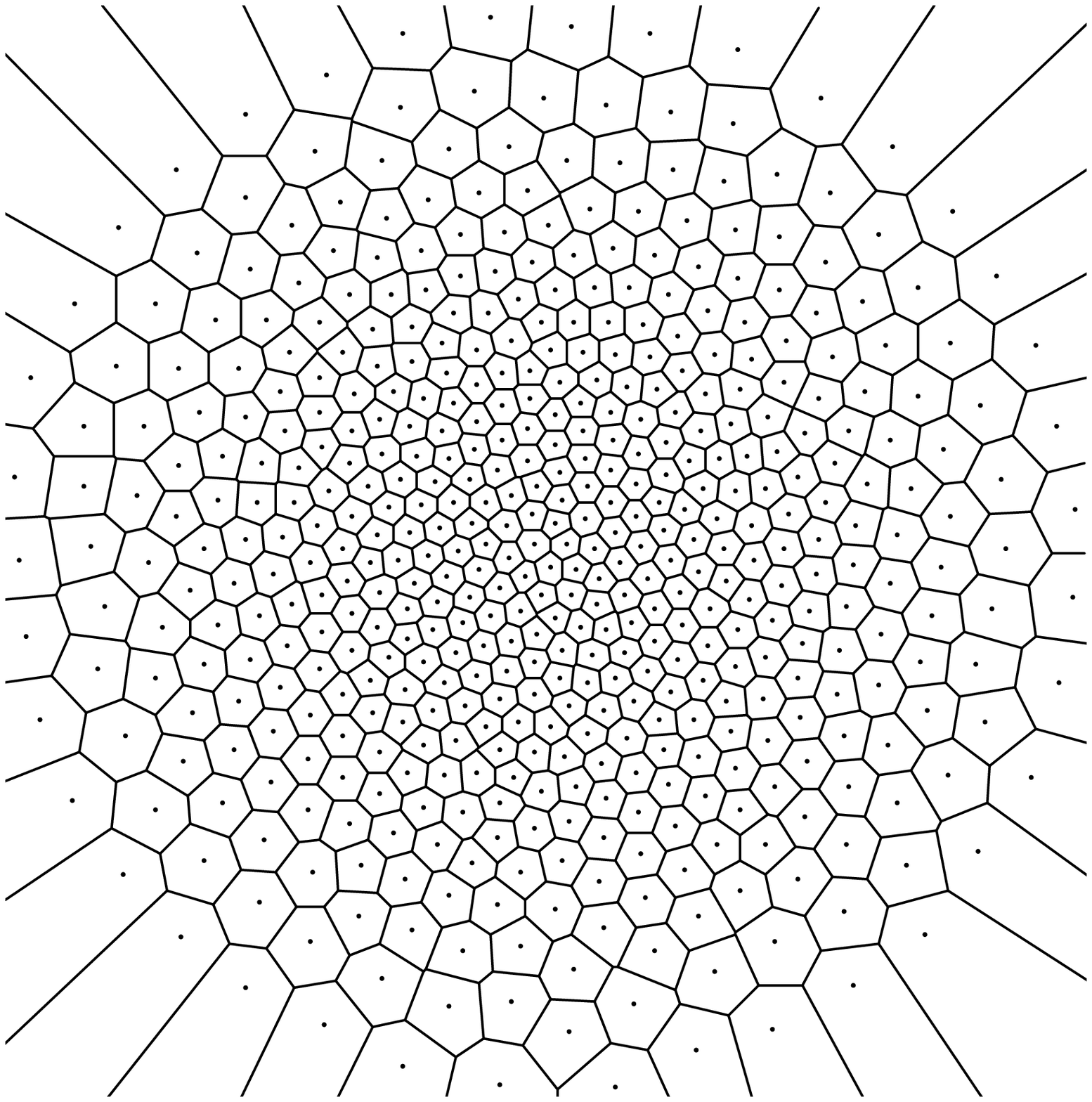}
\end{tabular}
\caption{\em Two $N$-quantizers (and their Voronoi diagram) related
to bi-variate normal distribution ${\cal N}(0;I_2)$ ($N=500$); which
one is the best?}
\end{figure}
\end{center}
It is convenient at this stage to make a correspondence between
quantizers of size at most $N$ and $N$-tuples of $H^N$: to any
$N$-tuple $x:=(x_1,\ldots,x_N)$ corresponds a quantizer $\Gamma:=\Gamma(x)=\{x_i,
\,i=1,\, \ldots,N\}$ (of size at most $N$). One introduces the
quadratic distortion, denoted $D^X_{_N}$, defined on $H^N$ as a (symmetric)
function  by
\begin{eqnarray*}
D^X_{_N} &:& H^N \longrightarrow \R_+\\
(x_1,\ldots,x_{_N})&&\!\!\!\longmapsto \E \left(\min_{1\le i\le N}
|X-x_i|_{_H} ^2 \right).
\end{eqnarray*}
Note that, combining~(\ref{eNGammaXH})  and the definition of the distortion,  shows that
$$
D^X_{_N}(x_1,\ldots,x_{_N})= \E \left(\min_{1\le i\le N}
|X-x_i|_{_H} ^2 \right) = \E \left( d(X,\Gamma(x))^2 \right)=
\|X-\widehat X^{\Gamma(x)}\|_{_2}^2
$$
so that,
\[
e_{_N}(X,H)= \inf_{(x_1,\ldots,x_{_N})\in H^N}
\sqrt{D^X_{_N}(x_1,\ldots,x_{_N})}.
\]

The following proposition shows the existence of an optimal   $N$-tuple $x^{(N,*)}\!\in
H^N$ such that $e_{_N}(X,H) =\sqrt{D^X_{_N}(x^{(N,*)})}$.  The corresponding optimal quantizer at level $N$ is denoted
$\Gamma^{(N,*)}:=\Gamma(x^{(N,*)})$. In finite dimension we refer to~\cite{POL} (1982) and in infinite dimension
to~\cite{CUMA} (1988) and
\cite{PAR} (1990); one may also  see~\cite{PAG0},~\cite{GRLU1} and~\cite{LUPA1}. For recent developments on existence and
pathwise regularity of optimal quantizer see~\cite{GRLUPA2}.

\begin{proposition}\label{optiquantizer}$(a)$ The function $D^X_{_N}$ is lower semi-continuous for
the  product weak  topology on $H^N$.

\smallskip
\noindent $(b)$ The function $D^X_{_N}$ reaches a minimum at a
$N$-tuple $x^{(N,*)}$ (so that $\Gamma^{(N,*)}$ is an {\em optimal quantizer at level
$N$}).

\smallskip
-- If ${\rm card}({\rm supp}(\P_{X}))\ge N$,  the quantizer has full
size $N$ ($i.e.$ ${\rm card}(\Gamma^{(N,*)})=N$) and
$e_{_N}(X,H)<e_{_{N-1}}(X,H)$.

\smallskip
-- If ${\rm card}({\rm supp}(\P_{X}))\le N$, $e_{_N}(X,H)=0$.

\smallskip
Furthermore $\displaystyle \lim_N e_{_N}(X,H)=0$.

\noindent $(c)$ Any optimal (Voronoi)
quantization at level $N$, $\widehat X^{\Gamma^{(N,*)}}$ satisfies
\begin{equation}\label{stationarity}
\widehat X^{\Gamma^{(N,*)}}= \E( X\,|\, \sigma(\widehat X^{\Gamma^{(N,*)}}))
\end{equation}
where   $\sigma(\widehat X^{\Gamma^{(N,*)}})$ denotes the $\sigma$-algebra generated by
$\widehat X^{\Gamma^{(N,*)}}$.

\smallskip
\noindent $(d)$ Any optimal (quadratic) quantization at level $N$ is a best least square ($i.e.$ $L^2(\P)$)
approximation of $X$ among all $H$-valued random variables taking at most $N$ values:
\[
e_{_N}(X,H) =
\|X-\widehat{X}^{\Gamma^{(N,*)}}\|_{_2}=
\min\{\|X-Y\|_{_2},\,Y:(\Omega,{\cal A})\to H,\, {\rm card}(Y(\Omega))\le N\}.
\]
\end{proposition}

\noindent {\bf Proof (sketch of)}: $(a)$ The claim follows from the
 l.s.c. of $\xi\mapsto |\xi|_{_H}$ for the weak topology and Fatou's Lemma.

\smallskip
\noindent $(b)$ One proceeds by  induction on  $N$. If $N=1$, the
optimal  $1$-quantizer is $x^{(N,*)} =\{ \E \,X \}$ and
$e_2(X,H)=\|X-\E \,X\|_{_2}$.

Assume now that an optimal quantizer $x^{(N,*)}=(x_1^{(N,*)},\ldots,x_{_N}^{(N,*)})$ does exist at
level $N$.

\smallskip
$\;$ -- If ${\rm card}({\rm supp}(\P))\le N$, then the $N+1$-tuple 
$(x^{(N,*)},x^{(N,*)}_{_N})$ (among other possibilities) is also optimal at level $N+1$ and
$e_{_{N+1}}(X,H)=e_{_{N}}(X,H)=0$.

\smallskip
$\;$  -- Otherwise, ${\rm card}({\rm supp}(\P))\ge N+1$, hence
$x^{(N,*)}$ has pairwise distinct components and there exists $\xi_{N+1}\!\in
{\rm supp}(\P_{_X})\setminus\{x_i^{(N,*)},\,i=1,\ldots,N\}\neq
\emptyset$.

\smallskip
Then, with  obvious notations,
\[
D^X_{_{N+1}}((x^{(N,*)},\xi_{_{N+1}}))<
D^X_{_{N}}(x^{(N,*)}).
\]
Then, the set $F_{N+1}:=\left\{x\!\in
H^{N+1}\,|\,D^X_{_{N+1}}(x) \le
D^X_{_{N+1}}((x^{(N,*)},\xi_{_{N+1}})) \right\}$ is non empty, weakly closed  since $D^X_{_{N+1}}$
is l.s.c.. Furthermore, it is bounded  in $H^{N+1}$. Otherwise there would exist a sequence $x_{(m)}\!\in H^{N+1}$ such
that $|x_{(m),i_m}|_{_H}= \max_i|x_{(m),i}|_{_H}\to +\infty$ as $m\to \infty$. Then, by Fatou's Lemma, one checks
that
\[
\liminf_{m\to\infty}D^X_{_{N+1}}(x_{(m)}) \ge D^X_{_{N}}(x^{(N,*)})
>D^X_{_{N+1}}((x^{(N,*)},\xi_{_{N+1}})).
\]
Consequently $F_{N+1}$  is  weakly compact and the  minimum of $D^X_{_{N+1}}$ on $F_{N+1}$  is
clearly its minimum over the whole space $H^{N+1}$. In particular
\[
e_{_{N+1}}(X,H)\le
D^X_{_{N+1}}((x^{(N,*)},\xi_{_{N+1}})) < e_{_N}(X,H).
\]
If ${\rm card}({\rm supp}(\P))= N+1$, set $x^{(N+1,*)}= {\rm
supp}(\P)$ (as sets) so that t
 $X= \widehat X^{\Gamma^{(N+1,*)}}$ which implies $e_{_{N+1}}(X,H)=0$.

 \smallskip To establish that $e_{_N}(X,H)$ goes to $0$, one
 considers an everywhere dense sequence $(z_k)_{k\ge 1}$ in the separable space $H$.
 Then, $d(\{z_1,\ldots,z_{_N}\},X(\omega))$ goes to $0$ as $N\to
 \infty$ for every $\omega\!\in\Omega$. Furthermore, $d(\{z_1,\ldots,z_{_N}\},X(\omega))^2\le |X(\omega)-z_1|_{_H}^2\!\in
 L^1(\P)$. One concludes by the Lebesgue  dominated convergence
 Theorem that $D^X_{_N}(z_1,\ldots,z_N)$ goes to $0$ as $N\to\infty$.

\smallskip
\noindent $(c)$ and $(d)$ Temporarily set $\widehat X^*:= \widehat
X^{\Gamma^{(N,*)}}$ for convenience. Let
$Y:(\Omega, {\cal A})\to H$ be a random vector taking at most $N$
values. Set $\Gamma:=Y(\Omega)$. Since    $\widehat X^{\Gamma}$ is a
Voronoi quantization of $X$ induced by $\Gamma$,
\[
|X-\widehat X^{\Gamma}|_{_H}=d(X,\Gamma)\le |X-Y|_{_H}
\]
so that
\[
\|X-\widehat X^{\Gamma}\|_{_2}\le \|X-Y\|_{_2}.
\]
On the other hand, the optimality of $\Gamma^{(N,*)}$ implies
\[
\|X-\widehat X^*\|_{_2}\le\|X-\widehat X^{\Gamma}\|_{_2}.
\]
Consequently
\[
\|X-\widehat X^*\|_{_2}\le\min\left\{\|X-Y\|_{_2},\;Y:(\Omega, {\cal
A})\to H,\;{\rm card }(Y(\Omega))\le N\right\}.
\]
The inequality holds as an equality since $\widehat X^*$ takes at most $N$ values. Furthermore,
considering random vectors of the form $Y= g(\widehat X)$ (which take at most as many values as the
size of $\Gamma^{(N,*)}$) shows, going back to the very definition of conditional expectation,  that $
\widehat X^*= \E(X\,|\,\widehat X^*)$ $\P$-$a.s.$ $\qquad_{\diamondsuit}$

\bigskip
Item~$(c)$ introduces a very important notion in (quadratic) quantization.
\begin{definition}  A quantizer  ${\Gamma}\subset H$ is {\em stationary} (or
{\em self-consistent}) if (there is a nearest neighbour projection such that $ \widehat{X}^{\Gamma} =
{\rm Proj}_{_\Gamma}(X)$
satisfying)
\begin{equation}\label{DefStatio}
  \widehat{X}^{\Gamma} = \E\left(X\,|\,
\widehat{X}^{\Gamma} \right).
\end{equation}
\end{definition}
Note  in particular that any stationary quantization satisfies $\E X=\E \widehat X^{\Gamma}$.

\medskip
As shown by Proposition~\ref{optiquantizer}$(c)$ any quadratic
optimal quantizer  at level $N$ is stationary. Usually, at least
when $d\ge2$, there are other stationary quantizers: indeed, the distortion function  
$D^X_{_N}$ is $|\,.\,|_{_H}$-differentiable at $N$-quantizers
$x\!\in H^N$ with pairwise distinct components and
\[
\nabla D^X_{_N} (x)= 2\left(\int_{C_i({x})}(x_i-\xi)
\P_{_X}\!(d\xi)\right)_{1\le i\le N}\hskip -0,55cm= \;2\left(\E (
\widehat X^{\Gamma(x)}-X)\mbox{\bf 1}_{\{\widehat{X}^{\Gamma(x)}=
 x_i\}}
 \right)_{1\le i\le N}.
\]
hence, any  critical points of $D^X_{_N}$ is a stationary  quantizer.

\bigskip
\noindent{\bf Remarks and comments.} $\bullet$ In fact (see~Theorem~4.2, p. 38, \cite{GRLU1}),  the Voronoi
partitions of $\Gamma^{(N,*)}$ always have  a  $\P_X$-negligible boundary so that (\ref{DefStatio}) holds for {\em any}
Voronoi diagram induced by $\Gamma$.

\smallskip
\noindent $\bullet$ The problem  of the uniqueness of
optimal quantizer (viewed as a set) is not mentioned in the above
proposition. In higher dimension, this essentially never occurs. In
one dimension, uniqueness of the optimal $N$-quantizer was first established
in~\cite{FLE}
with strictly $\log$-concave density function. This was successively
extended
in~\cite{KIE} and~\cite{TRU}
and lead to the following criterion (for
more general ``loss" functions than the square function):

\smallskip
\noindent If the distribution of $X$ is {\em absolutely continuous with a $\log$-concave
density function}, then, for every $N\ge 1$,  there exists  only one stationary quantizer of size $N$,
which turns out to be the optimal quantizer at level $N$.

\smallskip More recently,
a more geometric  approach to
uniqueness based on the  Mountain Pass Lemma  first developed in~\cite{LAPA} and then generalized
in~\cite{COH}) provided   a slight extension of the above criterion (in terms of loss functions).

\smallskip This $\log$-concavity assumption is satisfied by many families of probability distributions
like the uniform distribution on compact intervals, the normal
distributions, the gamma distributions. There are examples of
distributions with a non $\log$-concave density function having a
unique optimal quantizer for every $N\ge 1$ (see $e.g.$ the Pareto
distribution in~\cite{FOPA}). On the other hand simple examples of
scalar distributions having multiple optimal quantizers at a given
level can be found in~\cite{GRLU1}.

\smallskip
\noindent $\bullet$ A  stationary quantizer can be sub-optimal.
This will be emphasized in Section~\ref{OptiQuantB} for the Brownian motion
(but it is also true for finite dimensional Gaussian random vectors)
where some families of sub-optimal quantizers --~the {\em product
quantizers} designed from the Karhunen-Lo\`ve basis~-- are stationary quantizers.

\smallskip
\noindent $\bullet$ For the uniform distribution over an interval
$[a,b]$, there is a closed form for the optimal quantizer at level
$N$ given by
$\Gamma^{(N,*)}=\{a+(2k-1)\frac{b-a}{N},\,k=1,\ldots,N\}$. This
$N$-quantizer is optimal not only in the quadratic case but also for
any $L^r$-quantization (see a definition further on). In general there is no such closed form,
either in $1$ or higher dimension. However, in~\cite{FOPA}   some
semi-closed forms are obtained for several families of (scalar)
distributions including the exponential and the Pareto distributions: all the optimal quantizers can be expressed
using a single underlying sequence $(a_k)_{k\ge 1}$ defined by an induction
$a_{k+1}= F(a_k)$.

\smallskip
\noindent $\bullet$ In one dimension, as soon as the optimal quantizer at level $N$ is
unique (as a set or as an $N$-tuple with increasing components), it
is generally possible to compute it as the solution of the stationarity
equation~(\ref{stationarity}) either by a  zero search (Newton-Raphson gradient descent) or a fixed point
(like the specific Lloyd~I procedure, see~\cite{KIE2})  procedure.

\smallskip
\noindent $\bullet$ In higher dimension,   deterministic
 optimization methods become intractable and one uses stochastic
procedures to compute optimal quantizers. The main topic of this paper being functional quantization, we
postponed the short overview on these aspects to Section~\ref{OptiQuantB}, devoted to the optimal quantization
of the Brownian motion. But it is to be noticed  that all efficient optimization methods rely on the so-called {\em
splitting method}  which  increases progressively  the quantization level
$N$.  This method is directly inspired by the  
induction developed in the proof of claim~$(b)$ of Proposition~\ref{optiquantizer} since one designs the starting
value of the optimization procedure at size $N+1$ by ``merging" the
optimized $N$-quantizer obtained at level $N$ with one further point of $\R^d$, usually randomly sampled with
respect to an appropriate distribution (see~\cite{PAPR1} for a discussion).

\smallskip
\noindent  $\bullet$ As concerns functional quantization, $e.g.$ $H=L^2_{_T}$, there is a close connection
between the  regularity of optimal (or even stationary) quantizers  and that of $t\mapsto X_t$ form $[0,T]$ into
$L^2(\P)$. Furthermore, as concerns optimal quantizers of Gaussian processes, one shows (see~\cite{LUPA1}) that
they belong to the reproducing space of their covariance operator, $e.g.$ to the Cameron-Martin space
$H^1=\{\int_0^.\dot h_sds,\, \dot h \!\in L^2_{_T}\}$ when $X=W$.  Other properties of optimal quantization of
Gaussian processes are established in~\cite{LUPA1}.

\medskip
\noindent {\bf Extensions to the $L^r(\P)$-quantization of random variables.}
In this paper, we focus on the purely quadratic framework ($L^2_{_T}$ and $L^2(\P)$-norms), essentially because
it is a natural (and somewhat easier) framework for   the computation of optimized
grids for the Brownian motion and for some first applications (like the pricing of path-dependent options,
see section~\ref{pathdeppric}).  But a more general and natural framework is to consider the functional
quantization of   random vectors taking values in a separable Banach space $(E, |\,.\,|_{_E})$. Let
$X:(\Omega, {\cal A}, \P) \rightarrow (E,|\,\,\,|_{_E})$, such that $\E \,|X|_{_E}^r<+\infty$ for some $r\ge
1$ (the case $0<r<1$ can also be taken in consideration).

  The $N$-level $(L^r(\P),|\,.\,|_{_E})$-quantization problem for $X\!\in L_{_E}^r(\P)$ reads
\begin{eqnarray*}
e_{_{N,r}}(X,E)&:= &\inf\left\{\|X-\widehat X^{\Gamma}\|_{_r},\;
\Gamma \subset E,\; {\rm card}( \Gamma) \le N \right\}.
\end{eqnarray*}

The main examples for $(E, |\,.\,|_{_E})$ are the  non-Euclidean norms on  $\R^d$, the functional spaces
$L^p_{_T}(\mu):=L^p([0,T],\mu(dt))$, $1\le p\le \infty$, equipped with
its usual norm, $(E, |\,.\,|_{_E})=({\cal C}([0,T]),\|\,.\,\|_{\rm
sup})$, etc. As concerns, the existence of an optimal quantizer, it 
holds true for reflexive Banach spaces (see P\"arna (90)) and  $E=
L^1_{_T}$, but otherwise it may fail even when $N=1$
(see~\cite{GRLUPA2}). In finite dimension, the Euclidean feature is not
crucial (see~\cite{GRLU1}). In the functional setting, many results
originally obtained in a Hilbert setting have been extended to the Banach
setting either for existence or regularity results (see~\cite{GRLUPA2})
or for rates
see~\cite{Dereich},~\cite{Dereichetal},~\cite{LUPA2},~\cite{LUPA4}.

\section{Cubature formulae: conditional expectation and numerical integration}\label{Cubature}

Let $F:H\longrightarrow \R$ be a continuous functional (with respect to the
norm $|\,.\,|_{_H}$) and let $\Gamma\!\subset H$ be an
$N$-quantizer. It is natural  to
approximate $\E(F(X))$ by $\E(F(\widehat{X}^{\Gamma}))$. This
 quantity $\E(F(\widehat{X}^{\Gamma}))$   is simply the finite
 weighted sum
\[
\E\, (F(\widehat{X}^{\Gamma})) = \sum_{i=1}^N F(x_i) \P(\widehat
X^\Gamma=x_i).
\]
Numerical computation of $\E\, (F(\widehat{X}^{\Gamma}))$ is possible  as soon as $F(\xi)$ can be computed at any
$\xi\!\in H$ and the distribution $(\P(\widehat X=x_i))_{1\le i\le N}$ of $\widehat{X}^{\Gamma}$  is known. The induced
quantization error
$\|X-\widehat{X}^{\Gamma}\|_{_2}$ is used to control the error
(see below). These quantities related to the quantizer $\Gamma$ are also called   {\em companion parameters}.

Likewise, one can consider {\em a priori} the $\sigma(\widehat X^\Gamma)$-measurable random variable
$F(\widehat{X}^{\Gamma})$ as a good approximation of the conditional
expectation $\E(F(X)\,|\,\widehat{X}^{\Gamma})$.

\subsection{Lipschitz functionals}
Assume that the functional $F$ is Lipschitz continuous on $H$. Then
\[
 \left| \E(F(X)\,|\,\widehat{X}^{\Gamma})-
F(\widehat{X}^{\Gamma})\right| \le [F]_{_{\rm Lip}}
\E(|X-\widehat{X}^{\Gamma}|\,|\,\widehat{X}^{\Gamma})
\]
so that, for every real exponent $r\ge1$,
\[
\| \E(F(X)\,|\,\widehat{X}^{\Gamma})-
F(\widehat{X}^{\Gamma})\|_{_r} \le [F]_{_{\rm Lip}}
\|X-\widehat{X}^{\Gamma}\|_{_r}
\]
(where we applied conditional  Jensen inequality to the convex
function $u\mapsto u^r$). In particular, using that
$\E\,F(X)=\E(\E(F(X)\,|\,\widehat{X}^{\Gamma}))$, one derives (with 
$r=1$) that
\begin{eqnarray*}
 \left|\E\, F(X) -\E\,
F(\widehat{X}^{\Gamma})\right|& \le &\|
\E(F(X)\,|\,\widehat{X}^{\Gamma})-
F(\widehat{X}^{\Gamma})\|_{_1}\\
&\le& [F]_{_{\rm Lip}} \|X-\widehat{X}^{\Gamma}\|_{_1}.
\end{eqnarray*}
Finally, using the monotony of the $L^r(\P)$-norms as a function of $r$ yields
\begin{equation}\label{Cub1}
 \left|\E\,F(X) -\E\,
F(\widehat{X}^{\Gamma})\right| \le [F]_{_{\rm Lip}}
\|X-\widehat{X}^{\Gamma}\|_{_1}\le [F]_{_{\rm Lip}}\|X-
\widehat{X}^{\Gamma}\|_{_2}.
\end{equation}
In fact, considering the Lipschitz functional
$F(\xi):=d(\xi,\Gamma)$, shows that
\begin{equation}\label{L1carac}
\|X-\widehat{X}^{\Gamma}\|_{_1}=\sup_{[F]_{_{\rm Lip}}\le
1}\left|\E\, F(X) -\E\, F(\widehat{X}^{\Gamma})\right|.
\end{equation}
 The Lipschitz functionals making  up a characterizing family   for the weak convergence of probability measures on
$H$, one derives that, for any sequence of $N$-quantizers $\Gamma^N$ satisfying
$\|X-\widehat{X}^{\Gamma^N}\|_{_1}\to 0$ as $N\to \infty$,
\[
\sum_{1\le i\le N}
\P(\widehat{X}^{\Gamma^N}\!=x_i^N)\,\delta_{x^N_i}\stackrel{(H)}{\Longrightarrow}\P_{_X}
\]
where $\stackrel{(H)}{\Longrightarrow}$ denotes the weak convergence
of probability measures on $(H,|\,.\,|_H)$.

\subsection{Differentiable functionals  with Lipschitz differentials}  Assume now that  $F$ is
differentiable on $H$, with a Lipschitz continuous differential $DF$, and  that
the quantizer $\Gamma$ is {\em stationary} (see Equation~(\ref{DefStatio})).

A Taylor expansion yields
\begin{eqnarray*}
\left| F(X) - F(\widehat{X}^{\Gamma})-
DF(\widehat{X}^{\Gamma}).(X-\widehat{X}^{\Gamma}) \right|&\le&
[DF]_{_{\rm Lip}} |X-\widehat{X}^{\Gamma}|^2.
\end{eqnarray*}
Taking conditional expectation given $\widehat{X}^{\Gamma}$ yields
\begin{eqnarray*}
\left|\E (F(X)\,|\,\widehat{X}^{\Gamma}\!) \!-\! F(\widehat{X}^{\Gamma}\!)\!-\!
\E\left(DF(\widehat{X}^{\Gamma}\!).(X\!-\!\widehat{X}^{\Gamma}\!)\,|\,\widehat{X}^{\Gamma}\!\right)
\right|&\le& [DF]_{_{\rm Lip}}
\!\E(|X\!-\!\widehat{X}^{\Gamma}|^2 |\,\widehat{X}^{\Gamma}\!).
\end{eqnarray*}
Now, using that the random variable $DF(\widehat{X}^{\Gamma})$ is
$\sigma(\widehat{X}^{\Gamma})$-measurable, one has
\[
\E\left(DF(\widehat{X}^{\Gamma}).(X-\widehat{X}^{\Gamma})\right) =
\E
\left(DF(\widehat{X}^{\Gamma}).\E(X-\widehat{X}^{\Gamma}\,|\,\widehat{X}^{\Gamma})\right)=
0
\]
so that
\[
\left|\E (F(X)\,|\,\widehat{X}^{\Gamma}) - F(\widehat{X}^{\Gamma})
\right| \le  [DF]_{_{\rm Lip}}
\E\left(|X-\widehat{X}^{\Gamma}|^2\,|\,\widehat{X}^{\Gamma}\right).
\]
Then, for every real exponent $r\ge 1$,
\[
\left\|\E (F(X)\,|\,\widehat{X}^{\Gamma}) -
F(\widehat{X}^{\Gamma})\right\|_{_r}\le [DF]_{_{\rm Lip}}
\|X-\widehat{X}^{\Gamma}\|_{_{2r}}^2.
\]
In particular, when $r=1$, one derives like in the former setting
\begin{equation}\label{Cub2}
 \left|\E F(X) -\E
F(\widehat{X}^{\Gamma})\right| \le [DF]_{_{\rm Lip}}\|X-
\widehat{X}^{\Gamma}\|^{2}_{_2}.
\end{equation}
In fact, the above inequality holds provided $F$ is ${\cal C}^1$ with
Lipschitz differential on every Voronoi cell $C_i(\Gamma)$. A similar
characterization to~(\ref{L1carac}) based on these functionals could be established.

Some variant of these cubature formulae can be found in~\cite{PAPR1}
or~\cite{GRLUPA3} for functions or functionals $F$ having only
some local Lipschitz regularity.

\subsection{Quantized approximation of $\E(F(X)\,|\,Y)$}

Let $X$ and $Y$ be two $H$-valued random vector defined on the
same probability space $(\Omega,{\cal A},\P)$ and $F:H\to \R$ be  a
Borel functional. The natural idea is to approximate
$\E(F(X)\,|\,Y)$ by the quantized conditional expectation
$\E(F(\widehat X)\,|\,\widehat Y)$ where $\widehat X$ and $\widehat
Y$ are quantizations of $X$ and $Y$ respectively.

Let $\varphi_{_F}:H\to\R$ be a (Borel) version of the conditional
expectation $i.e.$ satisfying
\[
\E(F(X)\,|\,Y) =\varphi_{_F}(Y).
\]
Usually, no closed form is available for the function $\varphi_{_F}$ but
some regularity property can be established, especially in a
(Feller)  Markovian framework. Thus assume that both $F$ and
$\varphi_{_F}$ are Lipschitz continuous with Lipschitz coefficients $[F]_{\rm
Lip}$ and $[\varphi_{_F}]_{\rm Lip}$. Then
\[
\E(F(X)\,|\,Y)-\E(F(\widehat X)\,|\,\widehat Y)=\E(F(X)\,|\,
Y)-\E(F( X)\,|\,\widehat Y)+\E(F(X)-F(\widehat X)\,|\,\widehat Y).
\]
Hence, using that $\widehat Y$ is $\sigma(Y)$-measurable and that conditional expectation is an $L^2$-contraction,
\begin{eqnarray*}
\|\E(F(X)\,|\,
Y)-\E(F( X)\,|\,\widehat Y)\|_{_2}&=& 
\|\E(F(X)|Y)-\E(\E(F(\widehat X)| Y)|\widehat
Y)\|_{_2}\\ &\le&\|\varphi_{_{F}}\!(Y)-\E(F(  X)|\widehat Y)\|_{_2}
\\ &=&\|\varphi_{_{F}}\!(Y)-\E(\varphi_{_{F}}(Y)|\widehat Y)\|_{_2} 
\\
&\le&\|\varphi_{_{F}}\!(Y)- \varphi_{_{F}}(\widehat Y)\|_{_2}.
\end{eqnarray*}
The last inequality follows form the definition of conditional
expectation given $\widehat Y$ as the best quadratic approximation
among $\sigma(\widehat Y)$-measurable random variables.  On the other
hand, still using that $\E(\,.\,|\sigma(\widehat Y))$ is an $L^2$-contraction and this time that
$F$ is  Lipschitz continuous yields  
\[
\|\E(F(X)-F(\widehat X)\,|\,\widehat Y)\|_{_2}\le \|F(X)-F(\widehat
X)\|_{_2}\le [F]_{\rm Lip}\|X-\widehat X\|_{_2}.
\]
Finally, 
$$
 \|\E(F(X)\,|\,Y)-\E(F(\widehat X)\,|\,\widehat
Y)\|_{_2}\le [F]_{\rm Lip}\|X-\widehat X\|_{_2}+[\varphi_{_F}]_{\rm
Lip}\|Y- \widehat Y\|_{_2}.
$$

In the non-quadratic case the above inequality remains valid
provided $[\varphi_{_F}]_{\rm Lip}$  is replaced by $2[\varphi_{_F}]_{\rm
Lip}$.

\section{Vector quantization rate ($H=\R^d$)}\label{RateRd} The fact that
$e_{_{N}}(X,\R^d)$ is a non-increasing sequence that goes to $0$ as
$N$ goes to $\infty$ is a rather simple result
established in Proposition~\ref{optiquantizer}. Its rate of
convergence to $0$ is a much more challenging problem. An answer is
provided by the so-called Zador Theorem stated below.

This theorem was first stated  and established for distributions
with compact supports  by Zador (see~\cite{ZAD0,ZAD1}).
Then a first extension to general probability distributions on $\R^d$
is developed in~\cite{BUWI1}. The first mathematically rigorous proof  can be found in~\cite{GRLU1}, and relies on a
random quantization argument  (Pierce Lemma).
\begin{theorem}  $(a)$ {\sc Sharp rate}. Let $r>0$  and $X\!\in L^{r+\eta}(\P)$
for some $\eta>0$. Let $\P_{_X}(d\xi) = \varphi(\xi)\,d\xi
\stackrel{\perp}{+} \nu(d\xi)$  be the canonical decomposition of
the distribution of $X$ ($\nu$ and the Lebesgue measure are
singular). Then (if $\varphi\not \!\equiv 0$),
\begin{equation}\label{ZadorRate}
e_{_{N,r}}(X,\R^d)\sim \widetilde{J}_{r,d}\times\displaystyle
\left(\int_{\R^d}\varphi^{\frac{d}{d+r}}(u)\,du\right)^{\frac{1}{d}+\frac{1}{r}}\!\!\!\!\!\!\!\times
 N^{-\frac 1d} \quad \mbox{as}\quad
N \to +\infty.
\end{equation}
where $\widetilde{J}_{r,d}\!\in (0,\infty)$.

\medskip
\noindent $(b)$ {\sc Non asymptotic upper bound} (see $e.g.$~\cite{LUPA4}). Let $d\ge 1$. There exists
$C_{d,r,\eta}\!\in \!(0,\infty)$ such that, for every
$\R^d\!$-valued random vector $X$,
\[
\forall\,N\ge 1,\qquad e_{_{N,r}}(X,\R^d) \le
C_{d,r,\eta}\|X\|_{r+\eta} N^{-\frac 1d}.
\]
\end{theorem}

\noindent {\bf Remarks.}
 $\bullet$ The real constant $\widetilde{J}_{r,d}$ clearly corresponds to
 the case of the uniform distribution over the unit hypercube $[0,1]^d$ for which  the slightly more precise
statement holds
 \[
\lim_N N^{\frac 1d}e_{_{N,r}}(X, \R^d)=\inf_N N^{\frac 1d}e_{_{N,r}}(X, \R^d)=\widetilde
J_{r,d}.
 \]
 The proof is based on a self-similarity argument. The value of $\widetilde{J}_{r,d}$   depends
 on  the reference norm on $\R^d$. When $d=1$, elementary
 computations show that
 $\widetilde J_{r,1}= (r+1)^{-\frac 1r}/2$. When $d=2$, with the canonical Euclidean norm,
 one shows (see~\cite{NEW} for a proof, see also~\cite{GRLU1}) that
$\widetilde{J}_{2,d}=\sqrt{\frac{5}{18\sqrt{3}}}$.
 Its exact value is unknown for    $d\ge 3$ but, still for the canonical Euclidean
 norm, one has (see~\cite{GRLU1}) using some  random quantization arguments,
$$
\displaystyle \widetilde{J}_{2,d} \sim \sqrt{\frac{d}{2 \pi e}}\approx
\sqrt{\frac{d}{17,08}}\quad\mbox{ as }\quad d\to +\infty.
$$

\noindent $\bullet$ When  $\varphi\equiv 0$ the distribution of $X$ is purely singular. The   rate~(\ref{ZadorRate})
still holds in the sense that $\lim_N N^{\frac 1d}e_{_{r,N}}(X, \R^d)=0$.
 Consequently, this is not the right  asymptotics. The quantization problem
 for singular measures (like uniform distribution on fractal compact sets)
 has been extensively investigated by several authors, leading to
 the definition of a quantization dimension in connection with the rate
 of convergence of the quantization error on these sets. For more
 details we refer to~\cite{GRLU1, GRLU2} and the references therein.

\medskip
\noindent $\bullet$ A more naive way  to quantize the
uniform distribution on the unit
 hypercube is to proceed by {\em product quantization} $i.e.$ by
 quantizing the marginals of the uniform distribution. If $N= m^d$, $m\ge 1$, one easily proves that the best   quadratic
product quantizer (for the canonical Euclidean norm on $\R^d$) is the ``midpoint square grid"
\[
\displaystyle \Gamma^{sq,N} =\left(\frac{2i_1-1}{2m},\ldots,
\frac{2i_d-1}{2m}\right)_{1\le i_1,\ldots,i_d\le m}
\]
which induces a quadratic quantization error equal to
\[
\sqrt{\frac{d}{12}}\times N^{-\frac 1d}.
\]
Consequently, product quantizers are still {\em rate optimal} in
every dimension $d$. Moreover, note that the ratio of these two rates remains bounded as  $d\uparrow \infty$.

\section{Optimal quantization and  $QMC$}

The principle of Quasi-Monte Carlo method ($QMC$)  is  to approximate the integral  of a function $f:[0,1]^d\to \R$ 
with respect to the uniform distribution on
$[0,1]^d$, $i.e.$ $\displaystyle \int_{[0,1]^d}f \,d\lambda_d=\int_{[0,1]^d}f(\xi^1,\ldots,\xi^d)d\xi^1\cdots d\xi^d$ ($\lambda_d$ denotes the Lebesgue
measure on $[0,1]^d$), by the uniformly weighted  sum
\[
\frac 1N \sum_{k=1}^N f(x_k)
\] 
of values of $f$ at the points of a so-called low discrepancy $N$-tuple $(x_1,\ldots,x_{_N})$ (or set).
This  $N$-tuple can  the first $N$ terms of an infinite sequence.

If $f$ has finite variations denoted $V(f)$ --~either in the measure sense (see~\cite{BOLE,PAXI}) or in the Hardy and Krause sense (see~\cite{NIE}
p.19)~-- the Koksma-Hlawka inequality provides an  upper bound for the integration error induced by this method, namely
\[
\left|\frac 1N \sum_{k=1}^N f(x_k) -\int_{[0,1]^d}f
\,d\lambda_d\right|\le V(f)Disc^*_{_N}(x_1,\ldots,x_{_N})
\]
where  
$$
Disc^*_{_N}(x_1,\ldots,x_{_N}):= \sup_{y\in[0,1]^d}\left|\frac 1N \sum_{k=1}^N \mbox{\bf 1}_{\{x_k\in
[\![0,y]\!]\}}-\lambda_d([\![0,y]\!])\right|
$$
(with $[\![0,y]\!]\!=\!\prod_{k=1}^d[0,y^i]$, $y=(y^1,\ldots,y^d)\!\in[0,1]^d$).   The error modulus
$Disc^*_{_N}(x_1,\ldots,x_{_N})$ denotes {\em the discrepancy at the origin} of the $N$-tuple
$(x_1,\ldots,x_N)$.  For every $N\ge 1$, there exists  $[0,1]^d$-valued  $N$-tuples
$x^{(N)}$  such that 
\begin{equation}\label{lowboundDisc}
Disc^*_{_N}(x^{(N)})\le C_d \frac{(\log N)^{d-1}}{N},
\end{equation}
where $C_d\!\in(0,\infty)$ is a real constant only depending on $d$. This result can be proved using  the  so-called Hammersely procedure (see
$e.g.$~\cite{NIE}, p.31). When
$x^{(N)}=(x_1,\ldots,x_{_N})$ is made of the first $N$ terms of a  $[0,1]^d$-valued  sequence $(x_k)_{k\ge 1}$, then the above upper bound has be
replaced by  
$C'_d\frac{(\log N)^{d}}{N}$ ($C'_d\!\in (0,\infty)$). Such a sequence $x=(x_k)_{k\ge 1}$   is said to be a {\em sequence
with low discrepancy} (see~\cite{NIE} an the references therein for a comprehensive theoretical overview, but also~\cite{BOLE,PAXI} for examples
supported by numerical tests). When one only has
$Disc^*_{_N}(x_1,\ldots,x_{_N})\to 0$ as
$N\to\infty$, the sequence is said to be {\em uniformly distributed in} $[0,1]^d$. 

It is widely shared by
$QMC$ specialists that these rates are (in some sense) optimal although this remains a conjecture except when
$d=1$. To be precise what is known and what is conjectured is  the following:

\medskip
-- {\em Any} $[0,1]^d$-valued $N$-tuple  $x^{(N)}$ satisfies  $D^*_N(x^{(N)})\ge B_d
N^{-1}(\log N)^{\beta(d)}$ where $\beta(d) =\frac{d-1}{2}$ if $d\ge 2$ 
(see~\cite{ROT} and  also~\cite{NIE} and the references therein),
$\beta(1)=0$ and $B_d>0$ is a real constant only depending on $d$; the
conjecture is that
$\beta(d)=d-1$.

\medskip
\noindent 
-- {\em Any} $[0,1]^d$-valued sequence $(x_k)_{k\ge 1}$ satisfies $D^*_N(x^{(N)})\ge B_d
N^{-1}(\log N)^{\beta'(d)}$ for infinitely many $N$, where $\beta'(d)
=\frac{d}{2}$ if
$d\ge 2$  and $\beta'(1)=1$ and $B'_d>0$ is a real constant only depending on $d$; the conjecture is that $\beta(d)=d$. This follows from the result
for
$N$-tuple by the Hammersley procedure (see
$e.g.$~\cite{BOLE}).

\bigskip Furthermore, as concerns the use of Koksma-Hlawka inequality as an error bound for $QMC$ numerical integration, the different notions of
finite variation (which are closely connected) all become more and more restrictive --~and subsequently less and less ``natural" as a regularity
property of functions~--  when the dimension $d$ increases. Thus the
Lipschitz continuous function $f$ defined by 
$f(\xi^1,\xi^2,\xi^3):= (\xi^1+\xi^2+\xi^3)\wedge1$ has infinite
variation on $[0,1]^3$.

\medskip When applying Quasi-Monte Carlo approximation of  integrals with ``standard" continuous  functions on $[0,1]^d$, the best known error
bound, due to Proinov, is given by the following theorem.
\begin{theorem} (Proinov~\cite{PRO}) $(a)$ Assume $\R^d$ is equipped with the $\ell^\infty$-norm $|(u^1,\ldots,u^d)|_\infty:=\max_{1\le i\le
d}|u^i|$. Let $(x_1,\ldots,x_{_N})\!\in ([0,1]^d)^N$. For every  continuous function $f:[0,1]^d\to \R$, 
\[
\left|\int_{[0,1]^d} f(u)du-\frac 1N \sum_{k=1}^N f(x_k)\right|\le
K_d\,\omega_f((Disc^*_{_N}(x_1,\ldots,x_{_N}))^{\frac 1d})
\]
where $\omega_f(\delta):= \sup_{x,y\in [0,1]^d, |x-y|_\infty\le \delta}|f(x)-f(y)|$, $\delta\!\in(0,1)$,  
is the uniform continuity modulus of $f$ (with respect to the $\ell_\infty$-norm) and $C_d\!\in(0,\infty) $ is a universal constant only depending on $d$.  

\medskip
\noindent $(b)$ If $d=1$, $K_d=1$ and if $d\ge 2$, $K_d\!\in [1,4]$.
\end{theorem}

\noindent {\bf Remark.} Note that if $f$ is Lipschitz continuous, then $\omega_f(\delta)=[f]_{_{\rm Lip}}\delta$ where $[f]_{_{\rm Lip}}$ 
denotes the Lipschitz coefficient of $f$ (with respect to the $\ell_\infty$-norm).

\bigskip
First, this result emphasizes that low discrepancy sequences or sets do suffer from the {\em curse of dimensionality} when a $QMC$ approximation is
implemented  on functions having a ``natural" regularity like Lipschitz continuity. 

One also derives from this theorem an inequality between $(L^1(\P),\ell_\infty)$-quantization error of the uniform distribution $U([0,1]^d)$ and
the discrepancy at the origin of a $N$-tuple $(x_1,\ldots,x_{_N})$, namely
 \[
\|\,|U-\widehat U^{\{x_1,\ldots,x_{_N}\!\}}|_{\ell^\infty}\|_{_1}\le K_d (Disc^*_{_N}(x_1,\ldots,x_{_N}))^{\frac 1d}
\]  
since the function 
$\xi\mapsto \min_{1\le k\le N}|x_k-\xi|_{_\infty}$ is clearly $\ell_\infty$-Lipschitz continuous with Lipschitz coefficient $1$. The inequality
also follows from the characterization established in~(\ref{L1carac}) (which is clearly still true for non Euclidean norms). Then one may derive some
bounds for Euclidean norms (and in fact any norms) on $\R^d$ (probably not sharp in terms of constant) since all the norms  are strongly equivalent.
However the  bounds for optimal quantization error derived from Zador's Theorem ($O(N^{-\frac
1d}))$ and those for low discrepancy sets (see~(\ref{lowboundDisc}))
suggest that overall, optimal quantization provides lower error bounds
for numerical integration of Lipschitz functions than low discrepancy
sets, at least for  for generic values of
$N$. (However, standard computations show that for midpoint square grids
(with
$N=m^d$ points) both   quantization  errors and discrepancy    behave like $\frac 1m= N^{-\frac 1d}$).

\section{Optimal quadratic  functional quantization of Gaussian processes}\label{OptiQuantB}

Optimal quadratic functional quantization of Gaussian processes is closely
related to their so-called Karhunen-Lo\`eve expansion which
can be seen in some sense as some infinite dimensional  Principal Component Analysis ($PCA$) of a
(Gaussian) process.  Before stating a general result for Gaussian
processes, we start by the standard Brownian motion: it is the most important example in view of
(numerical) applications and for this process, everything can be made explicit.

\subsection{Brownian motion}
One considers the Hilbert space
$H= L^2_{_T}:=L^2([0,T],dt)$, $(f|g)_{_2} = \displaystyle \int_0^T\!\!f(t)g(t)dt$,
$|f|_{L^2_{_T}}= \sqrt{(f|f)_{_2}}$. The covariance operator $ C_{_W}$ of the Brownian motion $W=(W_t)_{t\in [0,T]}$ is defined on $L^2_{_T}$ by
\[
C_{_W}(f):= \E\left((f,W)_{_2}W \right)= \left(t\mapsto \int_0^T (s\wedge t)f(s) ds\right).
\]
\noindent  It is a symmetric positive trace class operator which can  be diagonalized  in the so-called  Karhunen-Lo\`eve
($K$-$L$) orthonormal basis $(e^W_n)_{n\ge 1}$ of $L^2_{_T}$, with eigenvalues $(\lambda_n)_{n\ge 1}$,   given by
\[
 e^W_n(t) = \sqrt{\frac 2T}\sin\left(\pi(n-\frac 12)\frac{t}{T}\right),\quad \lambda_n = \left(\frac{T}{\pi(n-\frac
12)}\right)^2, \; n\ge 1.
\]
This classical result can be established as a simple exercise by solving the functional equation $C_{_W}(f)=\lambda f$.
In particular, one can expand $W$ itself on this basis so that
\begin{eqnarray*}
W& \stackrel{L^2_{_T}}{=}&\sum_{n\ge 1} (W|e^W_n)_{_2}\, e^W_n.
\end{eqnarray*}
Now, the orthonormality of the ($K$-$L$) basis implies, using Fubini's Theroem,
\[
\E((W|e^W_k)_{_2} (W|e^W_\ell)_{_2}) = (e^W_k|C_{_W}(e^W_\ell))_{_2}= \lambda_\ell\delta_{k\ell}
\]
where $\delta_{k\ell}$ denotes the Kronecker symbol. Hence the Gaussian  sequence $((W|e^W_n)_{_2})_{n\ge 1}$ is pairwise non-correlated
which implies that these random variables are independent. The above identity also implies that ${\rm Var}((W|e^W_n)_{_2})=\lambda_n$.
Finally this shows that
\begin{equation}\label{KLexpansion}
W  \stackrel{L^2_{_T}}{=}  \sum_{n\ge 1} \sqrt{\lambda_n}\,\xi_n\,e^W_n
\end{equation}
where $\xi_n := (W|e^W_n)_{_2}/\sqrt{\lambda_n}$, $n\ge1$, is an i.i.d. sequence of ${\cal N}(0;1)$-distributed random variables.
Furthermore, this $K$-$L$ expansion converges in a much stronger sense since $\sup_{t\in [0,T]}|W_t-\sum_{k=1}^n \sqrt{\lambda_k}\xi_k
e^W_k(t)|\to 0$ $\P$-$a.s.$ and
\[
\|\sup_{[0,T]}|W_t-\sum_{1\le k\le n} \sqrt{\lambda_k}\xi_k e^W_k(t)|\|_{_2} = O\left(\sqrt{\log n/n}\right)
\]
(see~$e.g.$~\cite{LUPA3.5}). Similar results (with various rates) hold true for a wide class of Gaussian processes expanded on
``admissible" basis (see $e.g.$~\cite{LUPASpec}).

\begin{theorem}\label{QFW} (\cite{LUPA1} (2002) and~\cite{LUPA2} (2003)) Let
$\Gamma^N$, $N\ge 1$,  be a sequence of optimal $N$-quantizers for $W$.

\smallskip
\noindent $(a)$ For every $N\ge 1$, $ {\rm span}(\Gamma^N)
= {\rm span}\{e^W_1,\ldots,e^W_{d(N)}\}$ with $d(N)= \Omega(\log N)$. Furthermore $\widehat W^{\Gamma^N}$ and $W-\widehat
W^{\Gamma^N}$ are independent.

\medskip
\noindent $(b)$ $e_{_N}(W, L^2_{_T})=\displaystyle \|W-\widehat W^{\Gamma^N}\|_{_2}
 \sim \frac{T\sqrt{2}}{\pi}\frac{1}{\sqrt{\log N}}$ as $N\to\infty$.
\end{theorem}

\noindent {\bf Remark.} $\bullet$ The fact, confirmed by numerical experiments (see Section~\ref{OptiFBW}
Figure~\ref{dNegallogN}), that $d(N) \sim \log N$ holds as a conjecture.

\smallskip
\noindent $\bullet$ Denoting $\Pi_{_d}$ the orthogonal projection  on   ${\rm
span}\{e^W_1,\ldots,e^W_{d}\}$, one  derives from $(a)$  that $\widehat W^{\Gamma^N}=
\widehat{\Pi_{_{d(N)}}\!\!(W)}^{\Gamma_N}$ (optimal quantization at level $N$) and
\begin{eqnarray*}
 \|W-\widehat W^{\Gamma^N}\|^2_{_2}&=&\|\Pi_{d(N)}(W)- \widehat{\Pi_{_{d(N)}}\!\!(W)}^{\Gamma_N}\|^2_{_2} +
\|W-\Pi_{d(N)}(W)\|^2_{_2}\\ &=&  e_{N}\left(Z_{d(N)}, \R^{d(N)}\right)^2+\sum_{n\ge d(N)+1} \lambda_n
\end{eqnarray*}
where $ Z_{d(N)}\stackrel{d}{=}\Pi_{d(N)}(W) \sim \displaystyle
\bigotimes_{k=1}^{d(N)}  {\cal N}(0;\lambda_k)$.
\subsection{Centered Gaussian processes}

The above Theorem~\ref{QFW} devoted to the standard Brownian motion is a particular case of a more general theorem which holds for a wide
class of Gaussian processes

\begin{theorem}  (\cite{LUPA1} (2002) and  \cite{LUPA2} (2004)) Let $X=(X_t)_{t\in [0,T]}$ be a Gaussian process
with
$K$-$L$ eigensystem $(\lambda^X_n,e^X_n)_{n\ge
1}$ (with $\lambda_1\ge \lambda_2\ge\dots$ is non-increasing). Let $\Gamma^N$, $N\ge 1$,  be
a sequence of quadratic optimal $N$-quantizers for
$X$. Assume
\[
\lambda^X_n \sim \frac{\kappa}{n^b}\quad \mbox{ as } n\to \infty
\qquad (b>1).
\]

\noindent $(a)$  ${\rm span}(\Gamma^N)= {\rm span}\{e^X_1,\ldots,e^X_{d^X\!(N)}\}$ and $d^X(N)=
\Omega(\log N)$.

\noindent $(b)$ $e_{_N}(X, L^2_{_T})=\displaystyle \|X-\widehat
X^{\Gamma^N}\!\|_{_2} \sim
\sqrt{\kappa}\sqrt{b^b(b-1)^{-1}}\, (2\log N)^{-\frac{b-1}{2}}$.
\end{theorem}

\noindent{\bf Remarks.} $\bullet$  The above result admits an extension  to the case  $\displaystyle
\lambda^X_n\sim \varphi(n)$ as $n\to \infty$ with $\varphi$ regularly varying, index $-b\le -1$ (see~\cite{LUPA2}).
In~\cite{LUPA1}, upper or lower bounds are also established when
\[
(\lambda^X_n \le \varphi(n) ,\quad n\ge 1) \quad \mbox{ or }\quad (\lambda^X_n \ge \varphi(n),\quad n\ge 1).
\]
\noindent $\bullet$ The sharp asymptotics  $d^X(N)\sim
\frac 2b\log N$ holds as a conjecture.

\medskip \noindent {\em Applications to classical (centered) Gaussian processes.}
 
\smallskip
\noindent $\bullet$  Brownian bridge: $X_t := W_t-\frac{t}{T}W_T$, $t\!\in[0,T]$ and
$e^X_n(t) =\sqrt{2/T}\sin\left(\pi n\frac{t}{T}\right)$, $\lambda_n = \left(\frac{T}{\pi
n}\right)^2$, so that $e_{_N}(X, L^2_{_T})\sim T\frac{\sqrt{2}}{\pi}(\log N)^{-\frac 12}$.

\medskip
\noindent $\bullet$ Fractional Brownian motion  with Hurst constant $H\!\in (0,1)$
\[
e_{N}(W^{H}, L^2_{_T})\sim T^{H+\frac 12}c(H)(\log N)^{-H}
\]
where $c(H)= \left(\frac{\Gamma(2H) \sin(\pi H)(1+2H)}{\pi}\right)^{\frac 12}\!\!
\left(\frac{1+2H}{2\pi}\right)^H$
and $\Gamma(t)$ denotes the Gamma function at $t>0$.

\medskip
\noindent  $\bullet$ Some further explicit sharp rates can be derived from the  above
theorem for other classes of Gaussian stochastic processes (see~\cite{LUPA2}, 2004) like the fractional
Ornstein-Uhlenbeck processes, the  Gaussian diffusions,  a wide class Gaussian
stationary processes (the quantization rate is derived from the high frequency asymptotics of its spectral
density, assumed to be square integrable on the real line), for the $m$-folded integrated Brownian  motion,
the fractional Brownian sheet, etc.

\medskip
\noindent  $\bullet$ Of course some upper bounds can be derived for some even wider classes of processes,
based on the above first remark (see $e.g.$~\cite{LUPA1}, 2002).

\medskip
\noindent {\em Extensions to $r,p\neq 2$} When the processes have some self-similarity
properties, it is possible to obtain some sharp rates in the non purely quadratic case: this has
been done for fractional Brownian motion in~\cite{Dereichetal} using some quite different techniques
in which self-similarity properties plays there a crucial role. It leads to the following sharp
rates, for
$p\!\in[1,+\infty]$ and
$r\!\in(0,\infty)$
$$
e_{_{N,r}}(W^H, L^p_{_T})\sim  T^{H+\frac 12} c(r,H)(\log N)^{-H},\quad c(r,H)\!\in(0,+\infty).
$$
\subsection{Numerical optimization of quadratic functional quantization}\label{OptiFBW}

 Thanks to the scaling property of Brownian motion, one may focus on the normalized  case $T=1$. The numerical approach to
optimal quantization  of the Brownian motion is essentially  based on  Theorem~\ref{QFW} and the
remark that follows: indeed these results show that quadratic
optimal functional quantization of a centered Gaussian process
reduces to a finite dimensional optimal   quantization problem for a
Gaussian distribution with a diagonal covariance structure.   Namely
the optimization problem at level $N$ reads
\[
({\cal O}_N)\equiv \left\{\begin{array}{rcl}
\displaystyle e_{_N}(W,L^2_{_T})^2&:=& \displaystyle e_{_N}(Z_{d(N)},\R^{d(N)})^2  +\sum_{k\ge
d(N)+1} \lambda_k\\
\mbox{where }\quad  Z_{d(N)}&\stackrel{d}{=}  &\displaystyle\bigotimes_{k=1}^{d(N)} {\cal
N}(0,\lambda_k).
\end{array}\right.
\]
Moreover, if $\beta^N:=\{\beta^N_1,\ldots,\beta^N_{N}\}$ denotes an optimal $N$-quantizer of
$Z_{d(N)}$, then, the optimal $N$-quantizer
$\Gamma^N$ of $W$ reads $\Gamma^N=\{x^N_1,\ldots,x^N_N\}$ with
\begin{equation}\label{optizedquantW}
x^N_i(t) = \sum_{1\le \ell\le d(N)}  (\beta^N_i)^{\ell}e^W_\ell(t),\quad i=1,\ldots,N.
\end{equation}

The good news is that  $({\cal O}_N)$ is in fact a {\em finite
dimensional quantization optimization problem} for each $N\ge 1$.
The bad news is that the problem is somewhat ill conditioned  since the
decrease of the eigenvalues of $W$ is very steep for small values of
$n$:
$\lambda_1= 0.40528\dots$, $\lambda_2=0.04503\dots \approx
\lambda_1/10$. This is probably one reason for which former
attempts to produce good quantization of the Brownian motion first
focused on other kinds of quantizers like {\em scalar product
quantizers} (see~\cite{PAPR2} and Section~\ref{ProdFQ} below) or
$d$-dimensional block product quantizations (see~\cite{WIL}
and~\cite{LUPAWI}).

Optimization of the (quadratic) quantization of $\R^d$-valued random vector
has been extensively investigated since the early 1950's,
first in $1$-dimension, then in higher dimension when the cost of numerical
Monte Carlo simulation was drastically cut down
(see~\cite{GEGR}). Recent application of optimal vector quantization to numerics
turned out to be   much more  demanding  in terms of
accuracy. In that direction, one may cite~\cite{PAPR1}, \cite{MRBenH}
(mainly focused on numerical optimization of the quadratic
quantization of normal distributions). To apply the methods developed in
these papers, it is more convenient to rewrite our optimization
problem with respect to the standard $d$-dimensional distribution ${\cal N}(0;I_d)$
by simply considering the Euclidean norm derived
from the covariance matrix ${\rm Diag}(\lambda_1,\ldots, \lambda_{d(N)})$ $i.e.$

\[
({\cal O}_N)\Leftrightarrow\left\{\begin{array}{l}N\mbox{-optimal quantization of } \displaystyle
\bigotimes_{k=1}^{d(N)} {\cal
N}(0,1)\\
\mbox{for the covariance norm }
|(z_1,\ldots,z_{d(N)})|^2=
\sum_{k=1}^{d(N)} \lambda_k z^2_k.
\end{array}\right.
\]

The main point is of course that the dimension $d(N)$ is unknown.
However (see~Figure~\ref{dNegallogN}), one clearly verifies on   small values
of $N$ that the conjecture ($d(N)\sim \log N$) is most likely true. Then  for higher
values of $N$ one relies on it to shift from one dimension to
another following the rule $d(N) =d$,
$N\!\in\{e^{d},\ldots,e^{d+1}-1\}$.

\subsubsection{A toolbox for quantization optimization: a short overview}
Here is a short overview of stochastic optimization methods to compute
optimal or at least locally optimal quantizers in finite dimension. For
more details we refer to~\cite{PAPR1} and the references therein. Let $Z\stackrel{d}{=} {\cal
N}(0;I_d)$.

\medskip
\noindent{\em Competitive Learning Vector Quantization ($CLVQ$).}
This procedure is a  recursive stochastic approximation gradient descent
based on the integral representation of the gradient $\nabla
D^Z_N(x),\, x\!\in H^n$  (temporarily coming back to $N$-tuple notation)
of the distortion as the expectation of a {\em local gradient}
$i.e.$
$$
\forall\, x^N\!\in H^N,\quad\nabla D^Z_N(x^N) = \E (\nabla
D^Z_N(x^N,\zeta)),\;
\zeta_k\;\; i.i.d.,
\; \zeta_1\stackrel{d}{=} {\cal N}(0,I_d)
$$
so that, starting from  $x^N\!(0)\!\in (\R^d)^N$, one sets
\begin{eqnarray*}
\forall\,k\ge 0,\quad x^N\!(k+1)&=&
x^N\!(k)-\frac{c}{k+1}\nabla
D^Z_N(x^N\!(k) ,\zeta_{k+1})
\end{eqnarray*}
where $c\!\in(0,1]$ is a real constant to be tuned. As set, this looks quite formal but the
operating
$CLVQ$ procedure consists of two phases at each iteration:

\smallskip
$(i)$ {\em Competitive Phase:} Search of the nearest neighbor $x^N\!(k)_{i*(k+1)}$ of $\zeta_{k+1}$
among the components of
$x^N\!(k)_i$, $i=1,\ldots,N$  (using  a ``winning convention"  in case of conflict on the
boundary of the Voronoi cells).

\smallskip
$(ii)$ {\em Cooperative Phase:} One moves the winning component  toward $\zeta_{k+1}$ using a
dilatation $i.e.$
$x^N\!(k+1)_{i^*(k+1)}= {\rm Dilatation}_{\zeta_{k+1}, 1-\frac{c}{k+1}}( x^N\!(k)_{i^*(k+1)} )$.

\smallskip This procedure is useful for small or medium values of $N$. For an extensive study of this procedure, which turns out to be
singular in the world of recursive stochastic approximation algorithms, we
refer to~\cite{PAG1}. For general background on stochastic approximation,
we refer to~\cite{KUYI,BMP}.

\medskip
\noindent{\em The randomized ``Lloyd~I procedure".} This
is the randomization of the  stationarity based fixed point procedure since any optimal quantizer
satisfies~(\ref{DefStatio}):
\[
\widehat Z^{x^N\!(k+1)} = \E(Z\,|\,\widehat
Z^{x^N\!(k)}),\qquad
x^N\!(0)\subset \R^d.
\]
At every iteration the conditional expectation $\E(Z\,|\,\widehat
Z^{x^N\!(k)})$ is computed using a Monte Carlo simulation. For
more details about practical aspects of Lloyd~I procedure we refer
to~\cite{PAPR1}. In~\cite{MRBenH}, an approach  based on  genetic evolutionary algorithms is developed.

For both procedures, one may substitute a sequence of quasi-random numbers to the usual pseudo-random sequence.
This often speeds up the rate of convergence of the method,  although this can only be proved (see~\cite{LAPASA}) for a
very specific class of stochastic algorithm (to which $CLVQ$ does not belong).

\smallskip  The most important step to preserve the accuracy of the quantization as $N$ (and $d(N)$) increase is to
use the so-called {\em splitting method} which finds its origin in the proof
of the existence of an optimal $N$-quantizer: once the optimization
of a quantization grid of size $N$ is achieved, one specifies the
starting grid for the size $N+1$ or more generally $N+ \nu$, $\nu\ge 1$, by
merging the optimized grid of size $N$ resulting from the former
procedure with  $\nu$ points sampled independently from
the normal distribution with probability density proportional to
$\varphi^{\frac{d}{d+2}}$ where $\varphi$ denotes the p.d.f. of
${\cal N}(0;I_d)$. This rather  unexpected choice is motivated by
the  fact that this distribution provides  the
lowest    {\em in average}  random quantization error  (see~\cite{COH}).

\smallskip
 As a result, to be downloaded on the website~\cite{Website} devoted to quantization:

\medskip
\centerline{$\hbox{\tt www.quantize.maths-fi.com}$}

\medskip $\circ$   {\em Optimized  stationary codebooks for $W$}: in
practice, the $N$-quantizers $\beta^N$ of the distribution
$\otimes_{k=1}^{d(N)}{\cal N}(0;\lambda_k)$, $N\!=\!1$ up to
$10\,000$ ($d(N)$ runs from $1$ up to $9$).

\smallskip $\circ$   {\em Companion  parameters}:

\smallskip
$\quad$ -- distribution of $\widehat W^{\Gamma^N}$: $\P(\widehat
W^{\Gamma^N}\!=x^{N}_i)=\P(\widehat Z_{d(N)}^{\beta^{N}}\!= \beta^{N}_i)\; (
\leftarrow \!\!  \hbox{ in
 $\R^{d(N)}$})$.

$\quad$ -- The quadratic quantization error: $\|W-\widehat W^{\Gamma^N}\|_{_2}$.
\begin{figure}
\centering

 \begin{tabular}{cc}
\hskip -1.00 cm
\includegraphics[width=6cm,height = 1.5cm]{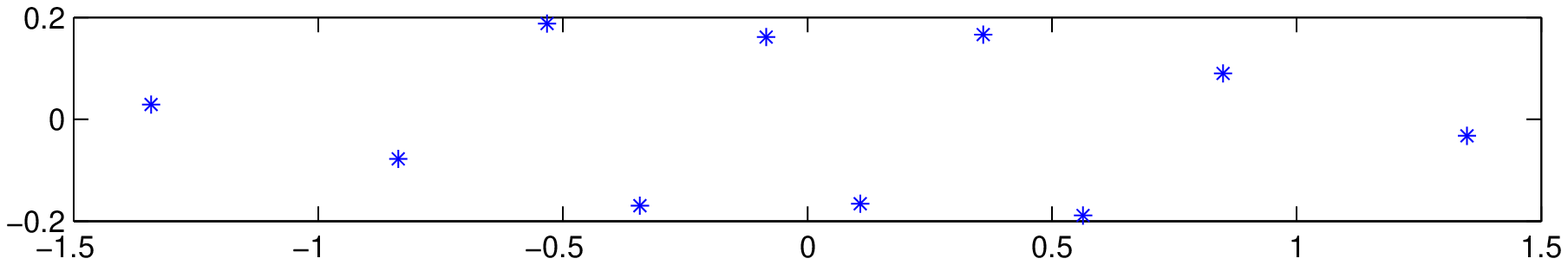}&\includegraphics[width=6cm,height=1.5cm]{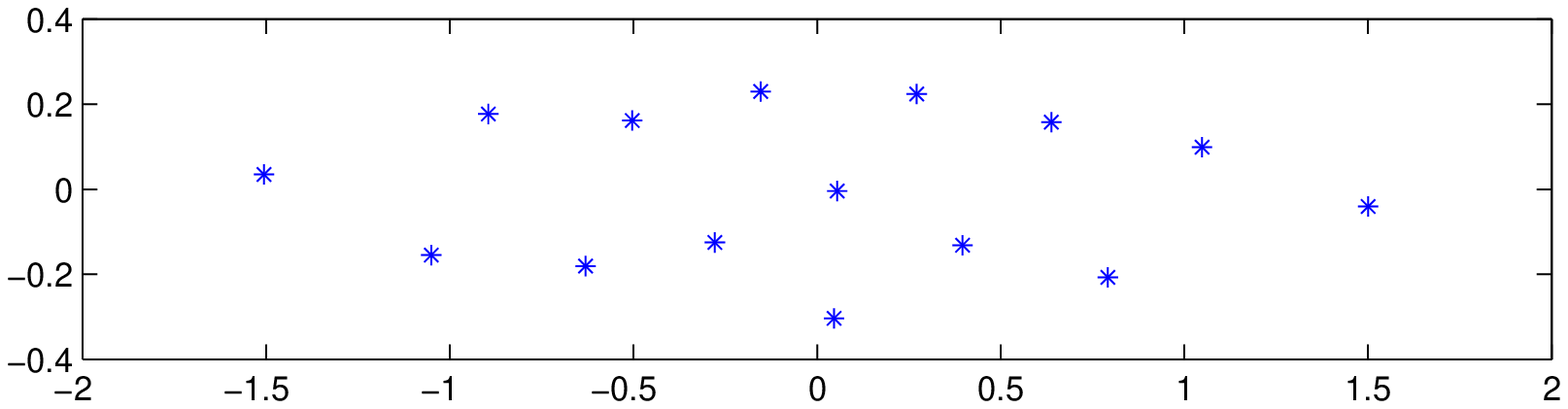}\\
\includegraphics[width=6cm,height = 4cm]{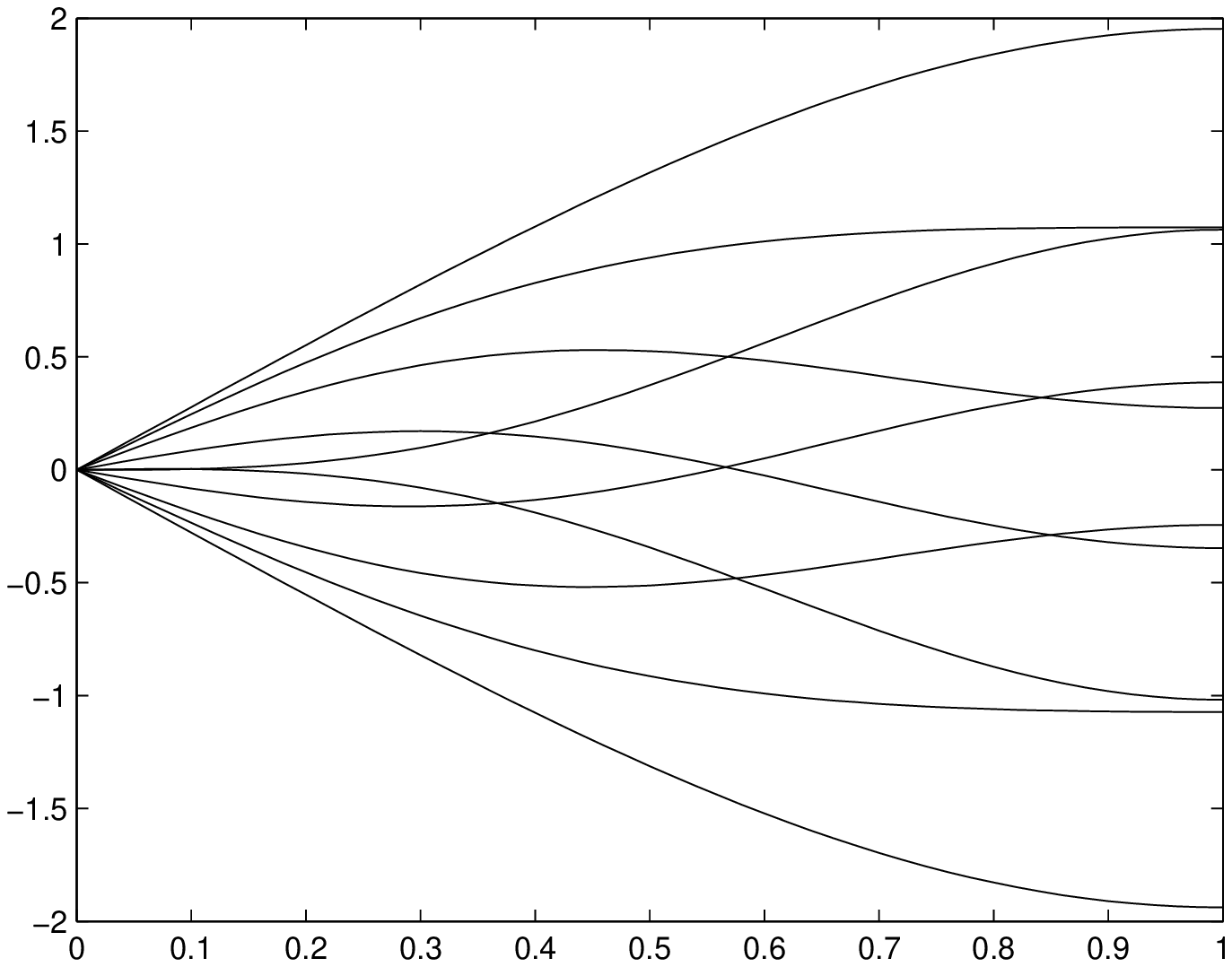}&\includegraphics[width=6cm,height= 4 cm]{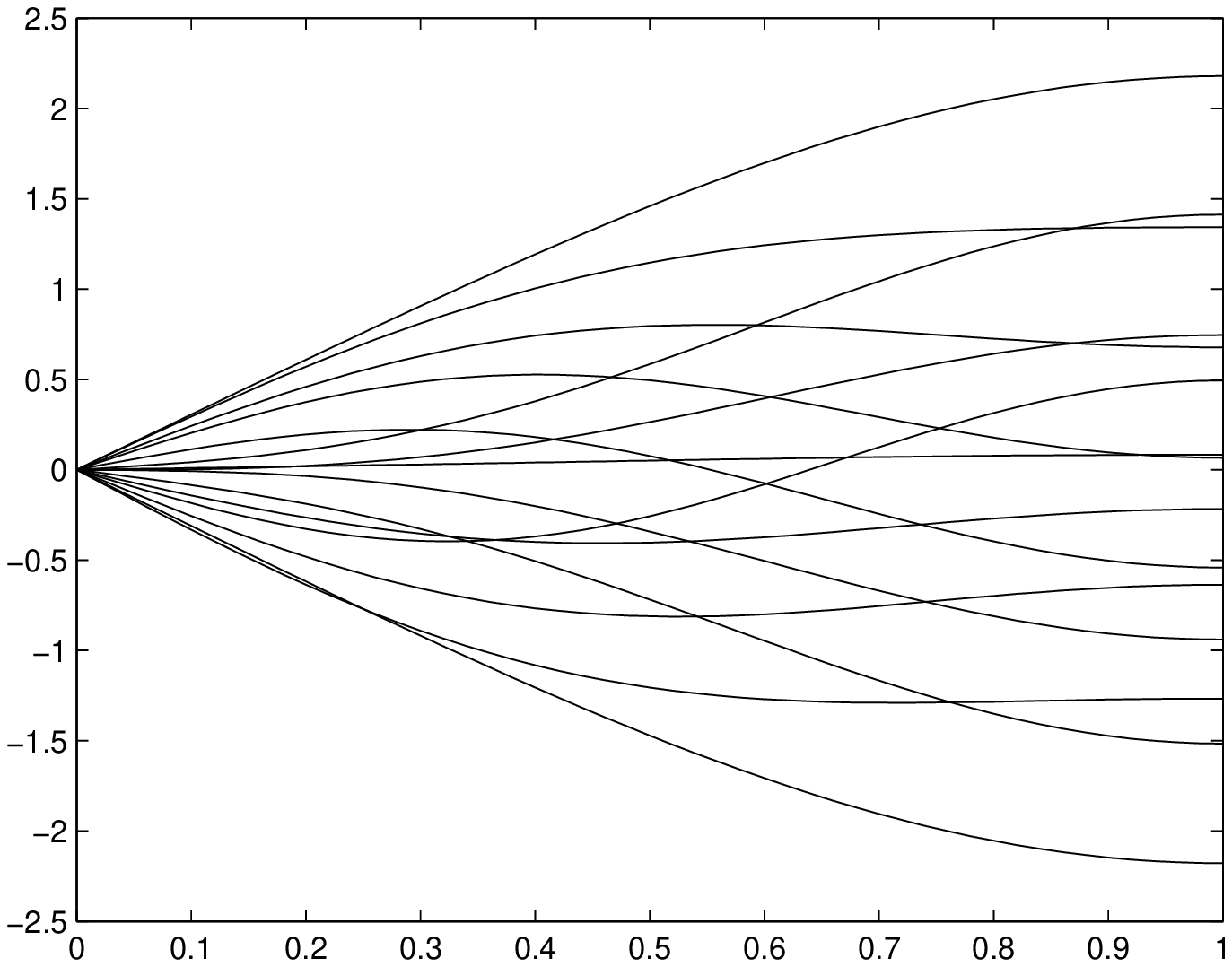}
\end{tabular}
\caption{\em Optimized functional quantization  of the Brownian
motion $W$ for $N=10,\,15$  ($d(N)=2$). Top: $\beta^{N}$ depicted in $\R^2$. Bottom:
 the optimized $N$-quantizer $\Gamma^{N}$.}
\end{figure}
 \begin{figure}
\centering
\begin{tabular}{ll}
\includegraphics[width=6cm,height = 4cm]{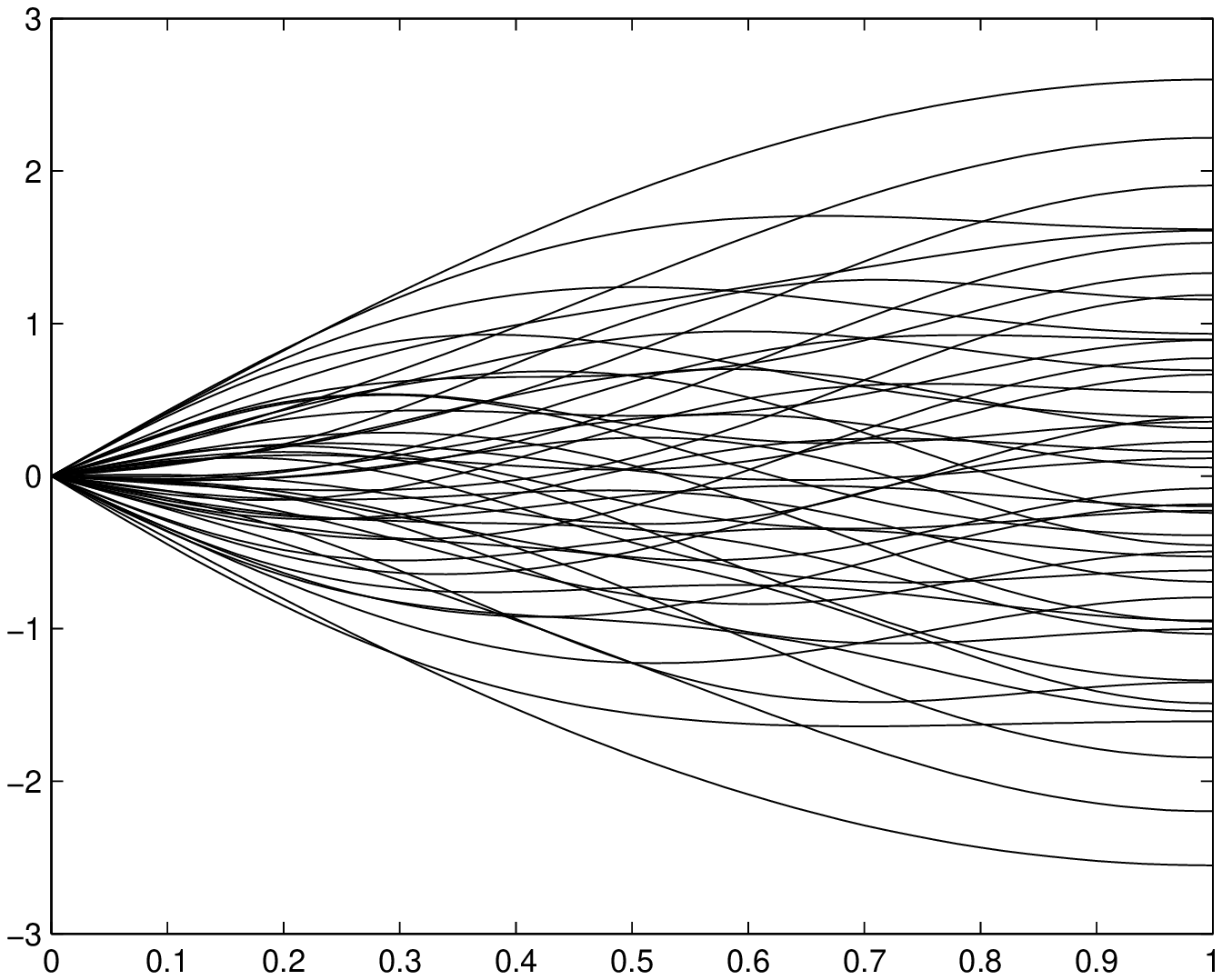}&\includegraphics[width=6cm,height = 4cm]{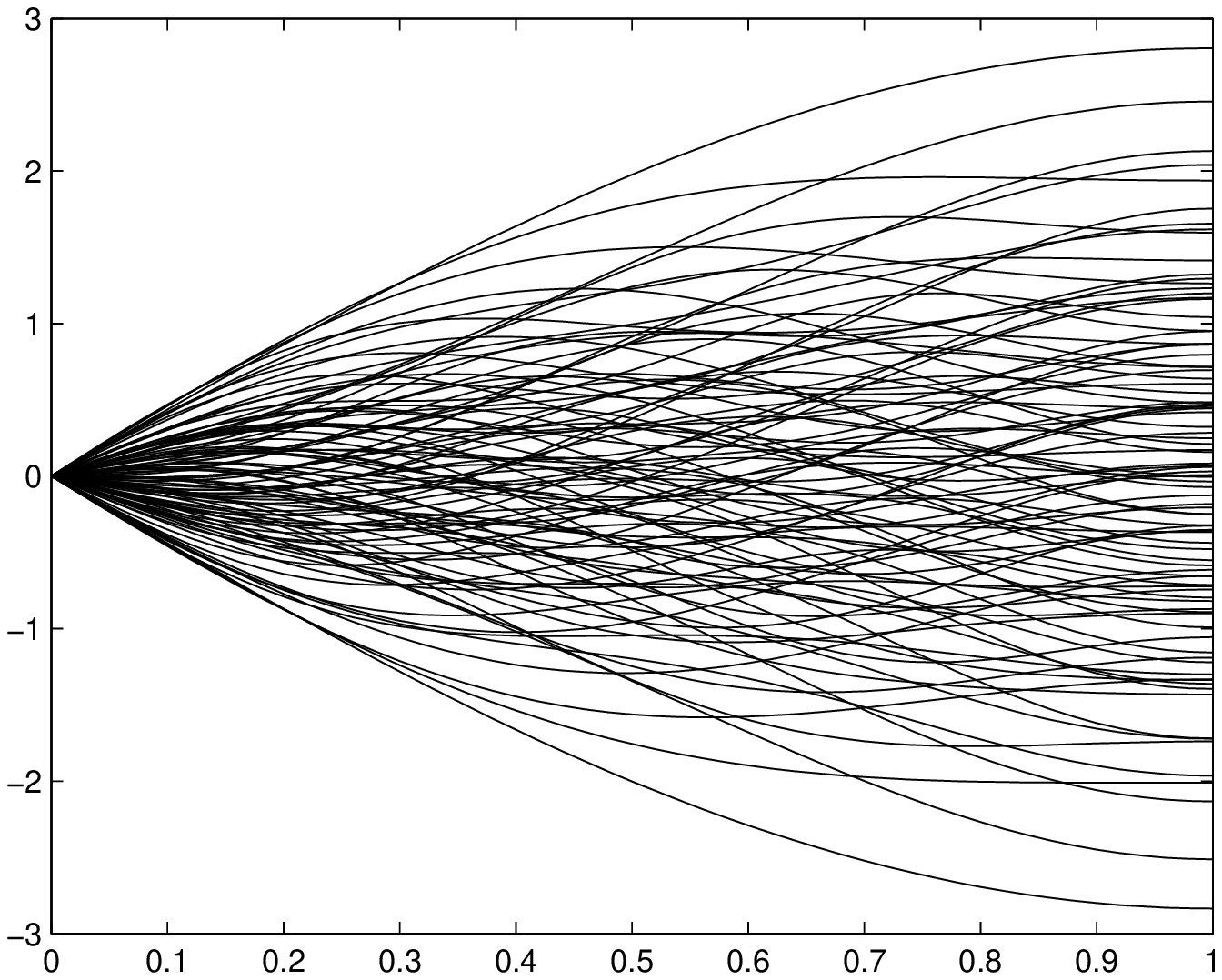}
\end{tabular}\caption{\em Optimized functional quantization   of the Brownian motion
$W$. The $N$-quantizers $\Gamma^{N}$. Left:  $N=48$ ($d(N)=3$). Right:  $N=96$, $d(96)=4$. }
\end{figure}
\begin{figure}
\centering
\begin{tabular}{c}
\hskip -0.75 cm  \includegraphics[width=13cm,height = 5cm]{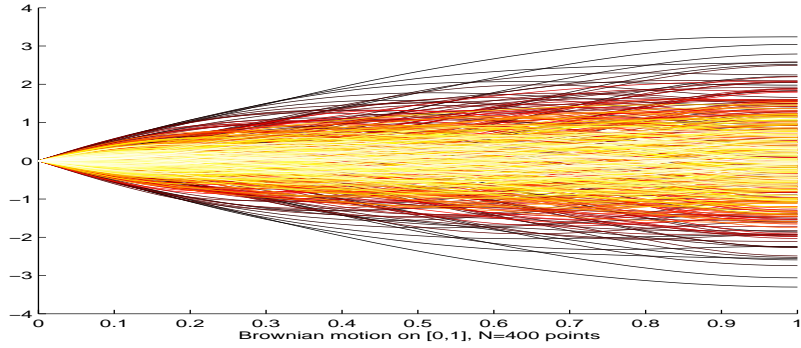}
\end{tabular}
\caption{\em Optimized $N$-quantizer $\Gamma^{N}$ of the
Brownian motion  $W$ with $N=400$. The grey level of the paths codes their  weights.}
\end{figure}
\begin{figure}\label{dNegallogN}
\centering
\begin{tabular}{c}
\includegraphics[width=8cm,height =5 cm ]{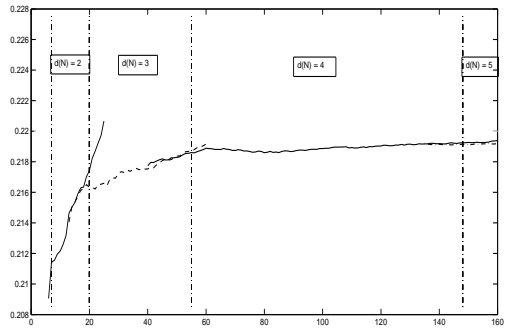}
\end{tabular}
\caption{\em  Optimal functional quantization of the Brownian motion.
$\displaystyle N\mapsto \log N \,(e_{_N}(W,L^2_{_T}))^2$,
$N\!\in\{6,\ldots,160\}$. Vertical dashed lines: critical dimensions for $d(N)$, $e^2 \approx 7$, $e^3\approx
20$, $e^4\approx 55$, $e^5\approx 148$.}
\end{figure}

\subsection{An alternative: product functional quantization}\label{ProdFQ}
Scalar Product functional quantization is a quantization method
which produces rate optimal  sub-optimal quantizers. They
were used $e.g.$ in~\cite{LUPA1} to provide exact rate (although not
sharp) for a very large class of processes. The first attempts to use functional quantization
for numerical computation with the Brownian motion   was achieved
with these quantizers (see~\cite{PAPR2}). We will see further on
their assets. What follows is presented for the Brownian motion but
would work for a large class of centered Gaussian processes.

 Let us consider again the expansion of $W$ in its $K$-$L$ basis :
\[
\;W  \stackrel{L^2_{_T}}{=}
\sum_{n\ge 1}
\sqrt{\lambda_n}\,\xi_n\,e^W_n
\]
where   $(\xi_n)_{n\ge 1}$ is an i.i.d. sequence ${\cal
N}(0;1)$-distributed random variables (keep in mind this convergence also holds $a.s.$ uniformly in $t\!\in[0,T]$). The idea is
simply to quantize these (normalized) random coordinates $\xi_n$: for every
$n\ge 1$, one considers  an optimal $N_n$-quantization of $\xi_n$,
denoted $\widehat{\xi}^{(N_n)}_n$ ($N_n\ge 1$). For $n> m$,
set $N_n=1$ and $\widehat{\xi}^{(N_n)}_n=0$ (which is the optimal
$1$-quantization). The integer  $m$ is  called the {\em length} of
the product quantization.  Then, one sets
\begin{eqnarray*}
\widehat W^{(N_1,\ldots,N_m,\,prod)}_t&\stackrel{}{:=}& \sum_{n\ge 1}
\sqrt{\lambda_n}\;\widehat{\xi}^{(N_n)}_n\,e^W_n(t) = \sum_{n= 1}^m
\sqrt{\lambda_n}\;\widehat{\xi}^{(N_n)}_n\,e^W_n(t).
\end{eqnarray*}
Such a quantizer takes $\prod_{n=1}^{m} N_n \le N$ values.

\medskip
 If  one denotes by
$\alpha^M=\{\alpha^M_1,\ldots,\alpha^M_M\}$ the (unique) optimal
quadratic $M$-quantizer of the ${\cal N}(0;1)$-distribution, the
underlying quantizer of the above quantization $\widehat
W^{(N_1,\ldots,N_m,\,prod)}$ can be expressed as follows (if one
introduces the appropriate multi-indexation): for every multi-index
$\underline{i}:=(i_1,\ldots,i_{m})\!\in
\prod_{n=1}^{m}\{1,\ldots,N_n\} $, set
\[
\displaystyle x^{(N)}_{\underline i}\!(t)\!:=\! \sum_{n=1}^{m}\!
\sqrt{\lambda_n}\,\alpha^{(N_n)}_{i_n} e^W_n\!(t)
\;\hbox{and}\;
\Gamma^{N_1,\ldots,N_m,prod}\!:=\!\left\{\!\displaystyle x^{(N)}_{\underline
i},\,\underline{i}\!\in \prod_{n=1}^{m}\{1,\ldots,N_n\}\!\! \right\}\!.
\]
Then the product quantization $\widehat W^{(N_1,\ldots,N_m,\,prod)}$
can be rewritten as
 \[
\widehat W^{(N_1,\ldots,N_m,\,prod)}_t=   \sum_{\underline{i}}
\mbox{\bf 1}_{\{W\in C_{\underline{i}}(\Gamma^{N_1,\ldots,N_m,prod})\}}
\displaystyle x^{(N)}_{\underline i}(t).
\]
where the Voronoi cell of $x^{(N)}_{\underline i}$ is given by
\[
C_{\underline{i}}(\Gamma^{N_1,\ldots,N_m,prod})  = \prod_{n=1}^m
(\alpha^{(N_n)}_{i_n-\frac 12},\alpha^{(N_n)}_{i_n+\frac 12})
\]
with $\alpha^{(M)}_{i\pm\frac
12}:=\frac{\alpha^{(M)}_{i}+\alpha^{(M)}_{i\pm 1}}{2}$,
$\alpha_0=-\infty$, $\alpha_{M+1}=+\infty$.

\subsubsection{Quantization rate by product quantizers}

It is clear that the optimal product quantizer is the solution to
the optimal integral bit allocation
\begin{equation}\label{OptiProdQuant}
\!\!\!\min\left\{\!\|W\!-\!\widehat
W^{(N_1,\ldots,N_m,\,prod)}\|_{_2},N_1,\ldots,N_m\ge 1, N_1\!\times
\!\cdots \!\times\! N_m\!\le\! N, m\!\ge \!1\!\right\}\!.\!
\end{equation}
Expanding $\|W-\widehat W^{(N_1,\ldots,N_m,\,prod)}\|_{_2}^2=
\||W-\widehat W^{(N_1,\ldots,N_m,\,prod)}|_{L^2_{_T}}\|_{_2}^2$ yields
\begin{eqnarray}
\label{identW1}\|W-\widehat W^{(N_1,\ldots,N_m,\,prod)}\|^2_{_2}&=&
\sum_{n\ge 1} \lambda_n \|\widehat
\xi_n^{(N_n)}- \xi_n\|_{_2}^2 \\
\label{identW2} &=& \sum_{n= 1} ^m\lambda_n (e_{_{N_n}}^2({\cal
N}(0;1),\R)-1)+\frac{T^2}{2}
\end{eqnarray}
\[
\mbox{since }\qquad \sum_{n\ge 1}\lambda_n =\E
\sum_{n\ge1}(W\,|\,e^W_n)_{_2}^2 =\E\int_0^TW_t^2dt=
\int_0^Tt\,dt=\frac{T^2}{2}.\quad\qquad\qquad
\]
\begin{theorem}(see~\cite{LUPA1}) For every $N\ge 1$, there exists an  optimal   scalar product quantizer
 of size at most $N$ (or at level $N$), denoted $\widehat W^{(N,\,prod)}$, of the
Brownian motion defined as the solution to the minimization
problem~(\ref{OptiProdQuant}). Furthermore these optimal product
quantizers make up a rate optimal sequence: there exists a real
constant $c_W>0$ such that
\[
\|W-\widehat W^{(N,\,prod)}\|_{_2}\le \frac{c_WT}{(\log N)^{\frac
12}}.
\]
\end{theorem}

\noindent {\bf Proof (sketch of).} By scaling one may assume without loss of generality that $T=1$. Combining~(\ref{identW1})
and Zador's Theorem shows
\begin{eqnarray*}
\|W-\widehat W^{(N_1,\ldots,N_m,\,prod)}\|^2_{_2} &\le& C
\left(\sum_{n= 1}^m \frac{1}{n^2N^2_n}\right)+\sum_{n\ge
m+1}\lambda_n\\
&\le& C' \left(\sum_{n= 1}^m \frac{1}{n^2N^2_n}+\frac 1m\right)
\end{eqnarray*}
with $\prod_n N_n \le N$. Setting $\displaystyle
m:=m(N)=\left[\log N\right]$ and $N_k= \left[\frac{(m!N)^{\frac
1m}}{k}\right]\ge 1$, $k=1,\ldots,m$, yields the announced
upper-bound. $\qquad_{\diamondsuit}$

\bigskip
\noindent {\bf Remarks.} $\bullet$ One can show that the length
$m(N)$ of the optimal quadratic product quantizer satisfies
\[
m(N) \sim \log N\qquad \mbox{ as }\qquad N\to +\infty.
\]

 \noindent $\bullet$  The most striking fact is that very few ingredients are
necessary to make the proof work as far as the quantization rate is
concerned. We only need   the basis of $L^2_{_T}$ on
which $W$ is expanded to be orthonormal {\em or} the random
coordinates to be orthogonal in $L^2(\P)$.  This robustness of the
proof has been used to obtain some upper bounds for very wide
classes of Gaussian processes by considering alternative orthonormal
basis of $L^2_{T}$ like the Haar basis for processes having
self-similarity properties (see~\cite{LUPA1}), or trigonometric
basis for stationary processes (see~\cite{LUPA1}). More recently,
combined with the non asymptotic Zador's Theorem, it was used to
provide some connections between mean regularity of stochastic
processes and quantization rate (see~Section~\ref{FQReg} and~\cite{LUPA4}).

 \smallskip
\noindent $\bullet$ Block quantizers combined with  large deviations estimates
 were used to provide the sharp rate obtained in Theorem~\ref{QFW} in~\cite{LUPA2}.

\smallskip
 \noindent $\bullet$ $d$-dimensional block quantization is also
 possible, possibly with varying  block size, providing a constructive approach to sharp
rate, see~\cite{WIL} and~\cite{LUPAWI}.

\smallskip
 \noindent $\bullet$ A similar approach can also provide some $L^r(\P)$-rates for product quantization with respect to the
$\sup$-norm over $[0,T]$, see~\cite{LUPA3.5}.

\subsubsection{How to use  product quantizers for numerical computations ?}

For numerics one can assume by a scaling argument  that $T=1$. To use product quantizers for numerics we need to have access to
the quantizers (or grid) at a given level $N$, their weights (and
the quantization error). All these quantities are available with
product quantizers. In fact the first attempts to use functional
quantization for numerics (path dependent option pricing) were
carried out with product quantizers (see~\cite{PAPR2}).

\medskip
 \noindent $\bullet$ The optimal product quantizers (denoted $\Gamma^{(N, prod)}$) at level $N$  are
 explicit,   given the optimal quantizers of the scalar normal
 distribution ${\cal N}(0;1)$. In fact the optimal allocation of the size $N_i$ of each
 marginal has been already achieved up to very high values of $N$. Some typical  optimal allocation (and the resulting
quadratic quantization error)  are reported in the table below.

\begin{center}\begin{tabular}{|c||c|c|c|}
\hline $N$&$N_{\rm rec}$&Quant. Error& Opti. Alloc.\\ \hline\hline
1&1&
0.7071&1\\ \hline 10&10&0.3138&5-2\\ \hline 100&96&0.2264&12-4-2  \\
\hline 1\,000&966&0.1881&23-7-3-2\\ \hline
 10\,000&9\,984 &0.1626& 26-8-4-3-2-2 \\   \hline
100\,000&   97\,920 &0.1461& 34 -- 10 -- 6 -- 4 -- 3 -- 2 -- 2 \\
\hline
\end{tabular}
\end{center}

\noindent $\bullet$ The weights $\P(\widehat W^{(N,\,prod)} =
x_{\underline i})$ are explicit too:
 the normalized coordinates $\xi_n$ of $W$ in its $K$-$L$ basis are
 independent, consequently
\begin{eqnarray*}
\P(\widehat W^{(N,\,prod)} = x_{\underline i})&=&\P(\widehat
\xi_n^{(N_n)}= \alpha^{(N_n)}_{i_n},\, n=1,\ldots,m(N) )\\
&= &\prod_{n=1}^{m(N)}\underbrace{\P(\widehat \xi_n^{(N_n)}=
\alpha^{(N_n)}_{i_n})}_{1D\,(tabulated)\,weights}.
\end{eqnarray*}

\noindent $\bullet$ Equation~(\ref{identW2}) shows that the
(squared)  quantization error of a product quantizer can be
straightforwardly computed as soon as one knows the eigenvalues and
the (squared)  quantization error of the normal distributions for
the $N_i$'s.

 The optimal allocations up to $N= 12\, 000$ can
be downloaded  on the website~\cite{Website}
as well as the necessary $1$-dimensional optimal quantizers (including the weights and the
quantization error) of the scalar normal
distribution (up to a size of $500$ which quite enough for this
purpose).

For numerical purpose we are also interested in the stationarity
property since such quantizers produce lower (weak) errors in  cubature formulas.

\begin{proposition}(see~\cite{PAPR2}) The product quantizers obtained from  the $K$-$L$ basis
are stationary quantizers (although sub-optimal).
\end{proposition}

\noindent {\bf Proof.} Firstly, note that
\begin{eqnarray*}
\widehat W^{N,prod}&=& \sum_{n\ge 1} \sqrt{\lambda_n}\, \widehat \xi^{(N_n)}_n e_n(t)
\end{eqnarray*}
so that $\sigma(\widehat W^{N,prod})=\sigma(\widehat \xi^{(N_k)}_k,\; k\ge 1)$. Consequently
\begin{eqnarray*}
 \E( W\,|\,\widehat W^{N,prod}) &=& \E(W\,|\, \sigma(\widehat \xi^{(N_k)}_k,\; k\ge 1))\\
  \E( W\,|\,\widehat W^{N,prod}) &=& \sum_{n\ge 1} \sqrt{\lambda_n}\,\E\left(\xi_n\,|\, \sigma(\widehat\xi^{(N_k)}_k,\; k\ge 1)\right)e^W_n\\
&\stackrel{i.i.d.}{=} & \sum_{n\ge 1} \sqrt{\lambda_n}\,\E\left( \xi_n\,|\, \widehat\xi^{(N_n)}_n\right)e^W_n\\
&=&  \sum_{n\ge 1} \sqrt{\lambda_n}\,\widehat \xi^{(N_n)}_ne^W_n
 \;=\; \widehat W.\qquad_\diamondsuit
\end{eqnarray*}

\noindent {\bf Remarks.} $\bullet$ This result is no longer true for
product quantizers based on other orthonormal basis.

\noindent $\bullet$ This shows the existence of non optimal   stationary quantizers.

\begin{figure}
\centering
\begin{tabular}{cc}
\includegraphics[angle=-90,width=6cm]{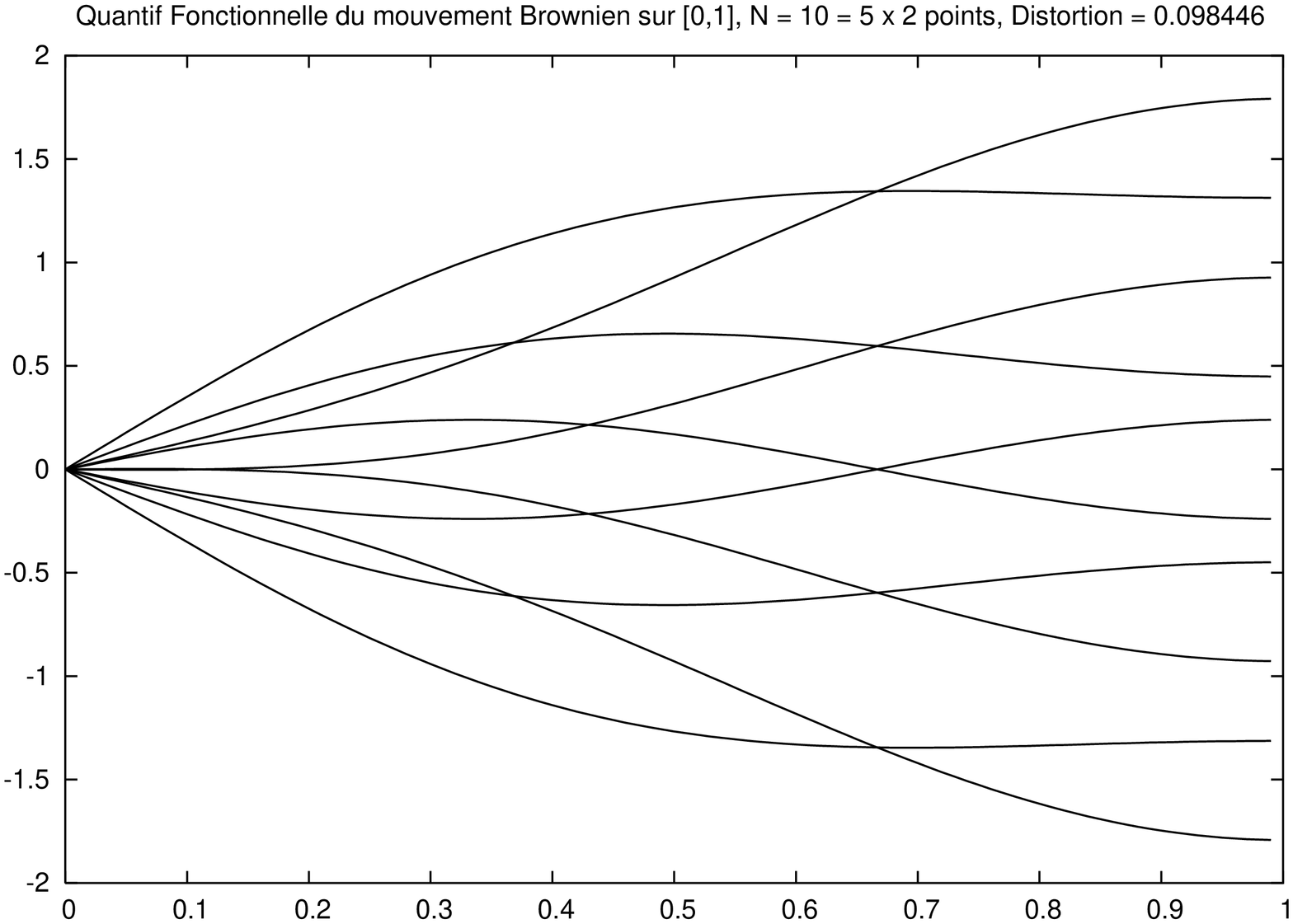}&
\includegraphics[angle=-90,width=6cm]{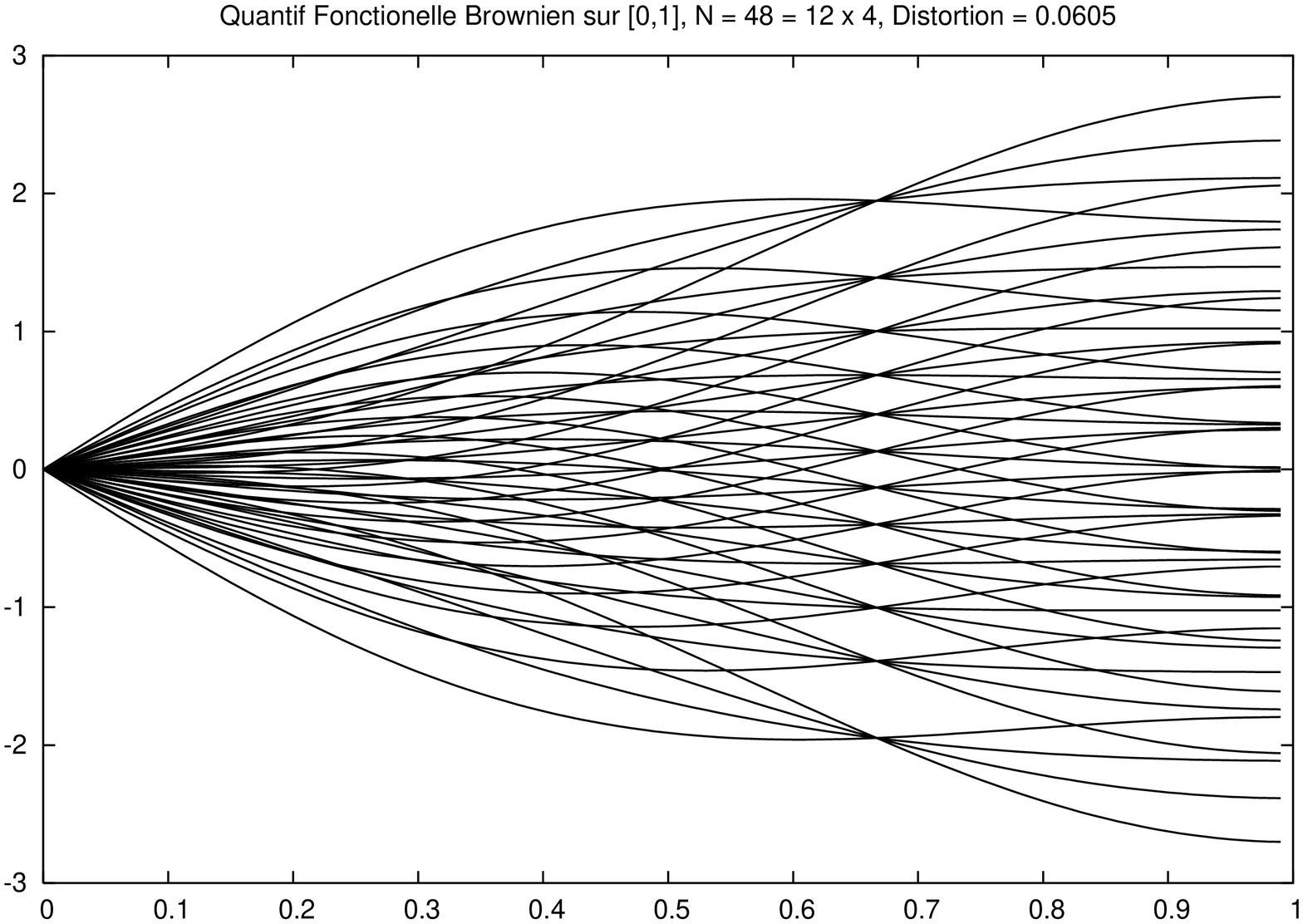}
\end{tabular}
\caption{\em  Product quantization of the Brownian motion: the $N_{\rm rec}$-quantizer $\Gamma^{(N,\,prod)}$.
$N=10$: $N_{\rm rec}=10$ and $N=50$: $N_{\rm rec}=12\times 4=48$.}
\label{fig:B10}
\end{figure}

\begin{figure}
\centering
\includegraphics[angle=-90,width=6cm]{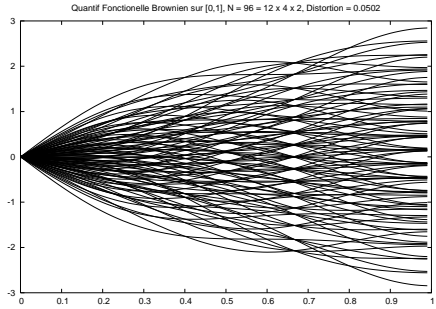}
\caption{\em Product quantization of the Brownian motion: the $N_{\rm rec}$-quantizer $\Gamma^{(N,\,prod)}$.
$N=100$: $N_{\rm rec}=12\times 4\times 2=96$.} \label{fig:B96}
\end{figure}

\subsection{Optimal $vs$  product quadratic functional quantization ($T=1$)}\label{OvsP}

\begin{figure}
\centering
\includegraphics[width=10.5cm,height=3.75cm]{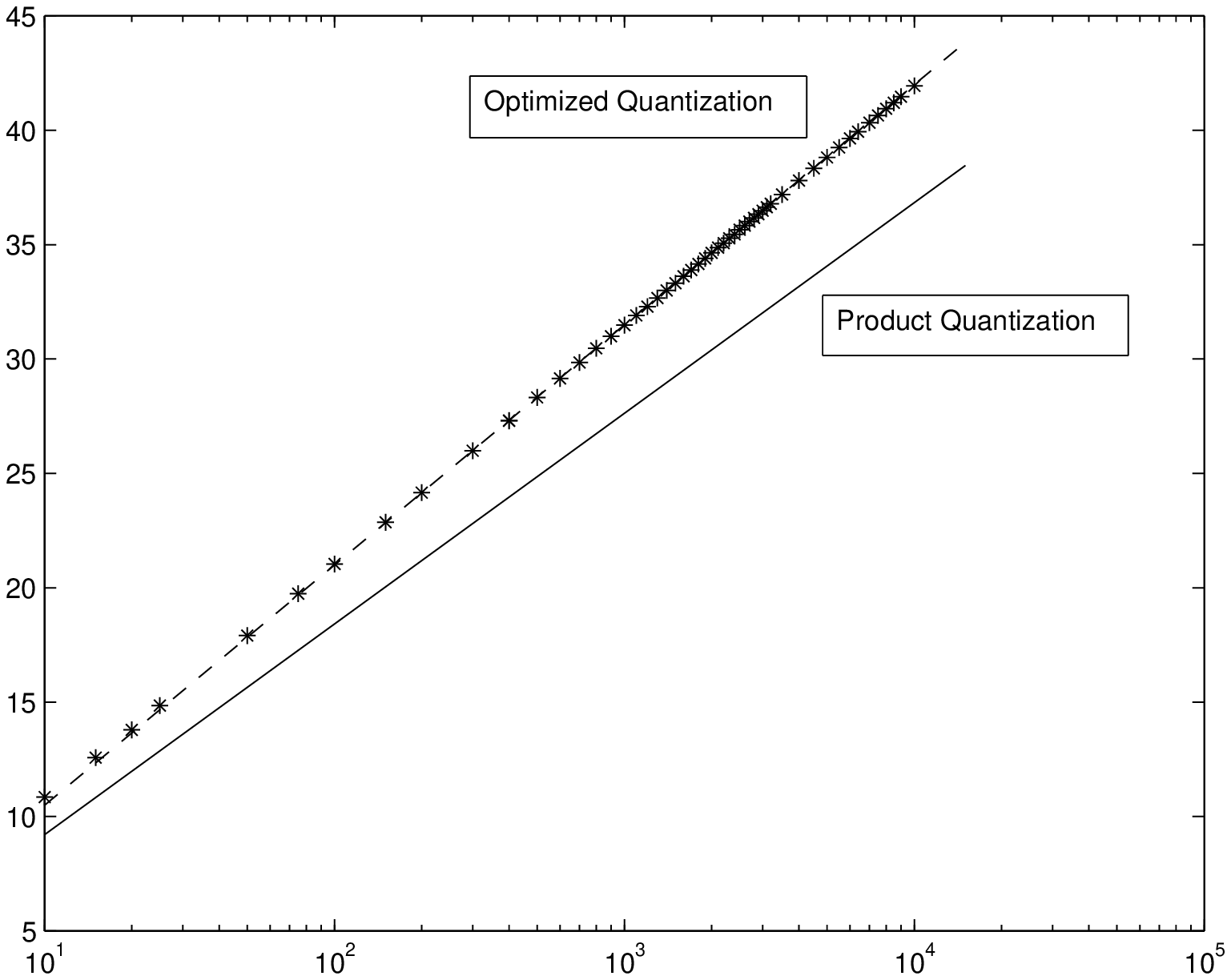}

\includegraphics[width=3.5 cm,height=9cm, angle= 270]{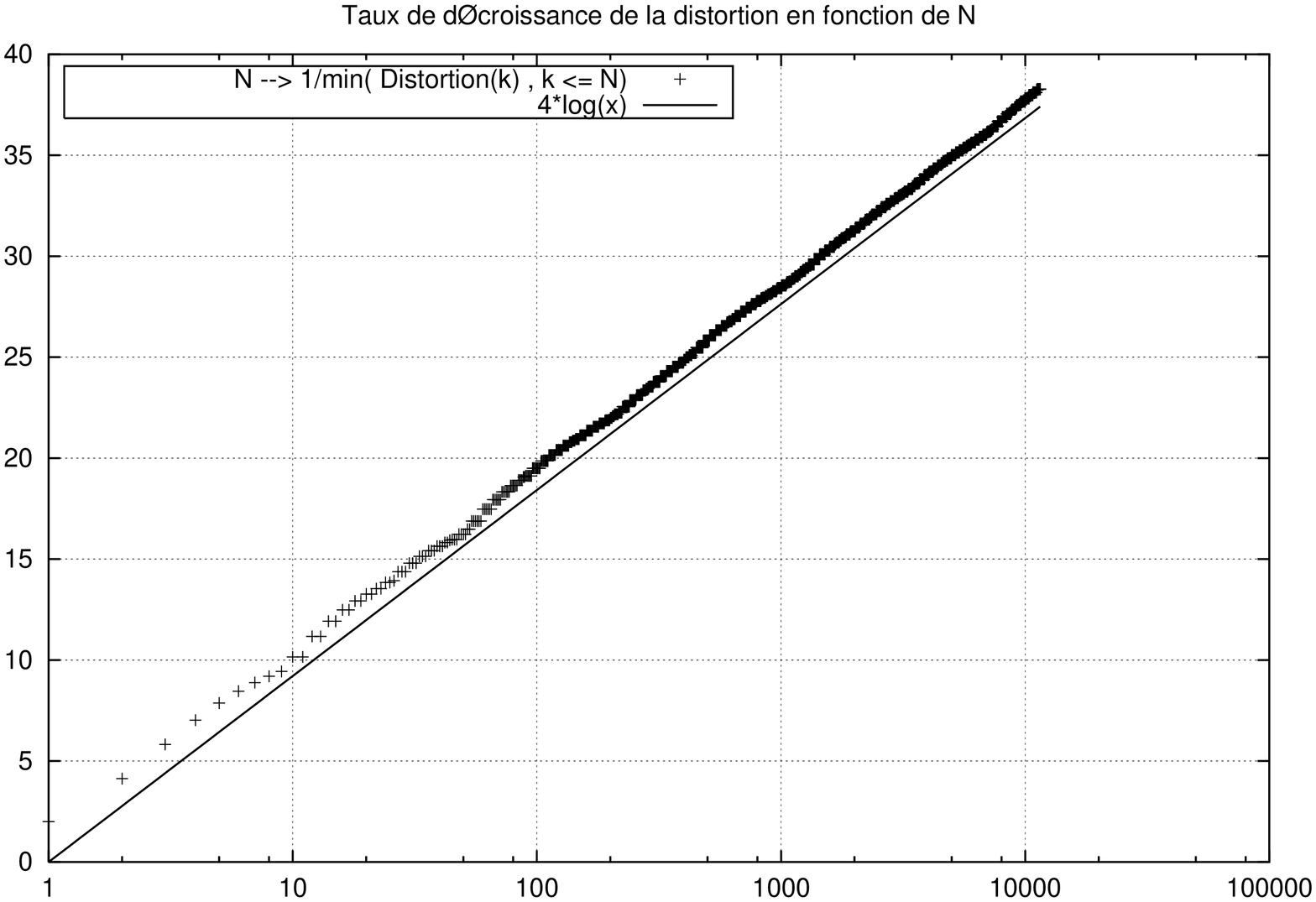}
\caption{\em    Numerical quantization rates. Top (Optimal
quantization).    Line~$\!+\!+\!+ $: $\log N\mapsto  (\|W-
\widehat W^{N}\|_{_2})^{-2}$. Dashed line: $\log N\mapsto \log
N/0.2194$. Solid line: $\log N\mapsto \log N/0.25$.  
Bottom (Product quantization). Line~$\!+\!+\!+$: $\log N \mapsto
\displaystyle (\min_{1\le k\le N}\|W-\widehat
W^{k,prod}\|^2_{_2})^{-1}$. Solid line: $\log N \mapsto
\log N/0.25$.}\label{fig10}
\end{figure}

$\circ$ {\sc (Numerical) Optimized Quantization:} By scaling, we can assume without loss of generality that $T=1$.
We carried out a huge optimization task in order to produce some {\em optimized} quantization
 grids for the Brownian motion by solving numerically   $({\cal O}_N)$ for $N=1$ up to $N= 10\,000$.
\[
e_{_N}(W,L^2_{_T})^2 \approx \frac{0.2195}{\log N}, \qquad
N=1,\ldots, 10\,000.
\]
This value (see~Figure~\ref{fig10}(left)) is significantly greater than the theoretical
(asymptotic) bound given by Theorem~\ref{QFW}  which is
$$
\lim_N \log Ne_{_N}(W,L^2_{_T})^2=\frac{2}{\pi^2}= 0.2026...
$$
Our guess, supported by our numerical experiments,  is that in fact $N\mapsto \log Ne_{_N}(W,L^2_{_T})^2$ is
possibly not monotonous but   unimodal.

\medskip
$\circ$   {\sc Optimal Product quantization}: as displayed  on Figure~\ref{fig10}(right), one has approximately
\[
\min\left\{\|\,|W-\widehat W|_{L^2_{_T}}\|^2_{_2},\, 1\le N_1\cdots
N_m\le N,\; m\ge 1\right\}=\|W- \widehat W^{(N,\,prod)}\|^2_{_2}\approx \frac{0.245}{\log N}
\]

$\circ$ {\sc Optimal $d$-dimensional block product quantization}: let us
 briefly mention this approach developed in~\cite{WIL} in which
 product quantization is achieved by quantizing some
 marginal blocks of size $1$, $2$ or $3$. By this approach, the
 corresponding constant is approximately $0. 23$, $i.e.$ roughly in
 between scalar product quantization and  optimized numeric quantization.

\medskip The conclusion, confirmed by our numerical experiments on option pricing (see~Section~\ref{pathdeppric}), is that

-- Optimal quantization is significantly more accurate
on numerical experiments but is much more demanding since it needs to keep off line or
at least to handle   large files (say 1 $GB$ for $N=10\,000$).

-- Both approaches  are included in the option pricer {\sc Premia} (MATHFI Project, Inria). An online  benchmark  is available on the website~\cite{Website}.

\vskip -2cm

\section{Constructive functional quantization of diffusions}
\subsection{Rate optimality for Scalar Brownian diffusions}

One considers on a probability space $(\Omega   , {\cal A}, \P)$ an homogenous Brownian diffusion
process:
$$
dX_t = b(X_t)dt +\vartheta(X_t)\,dW_t,\quad X_0=x_0\!\in\R,
$$
where $b$ and $\vartheta$ are continuous on $\R$ with at most linear growth ($i.e.$ $|b(x)|+|\sigma(x)|
\le C(1+|x|)$) so that at least a weak solution to the equation exists.

To devise a constructive way to quantize the diffusion $X$, it seems natural to start from  a
rate optimal quantization of the Brownian motion and to obtain some ``good" (but how good?)
quantizers for the diffusion by solving an appropriate $ODE$.  So let
$\Gamma^N=(w_1^N,\cdots, w_N^N)$, $N\ge 1$, be a sequence of
stationary rate optimal $N$-quantizers of $W$. One considers the following (non-coupled)
Integral Equations:
\begin{equation}\label{ODEd1}\displaystyle dx_i^{(N)}\!(t)= \left(b(x_i^{(N)}\!(t))
-\frac12 \vartheta\theta'(x_i^{(N)}\!(t))\right)dt
+\vartheta(t,x_i^{(N)}\!(t))\,d w_i^N(t).
\end{equation}

Set
\[
 \widetilde X^{ x^{(N)}}_t=\sum_{k=1}^N x_i^{(N)}\!(t)\mbox{\bf 1}_{\{\widehat W^{ \Gamma^N}=  w_i^N\}}.
\]
The process $\widetilde X^{ x^{(N)}}$ is a {\em non-Voronoi}   quantizer (since it is
defined using the   Voronoi diagram   of $W$). What is interesting is that it is
a {\em computable quantizer} (once the above integral equations have been solved) since the
weights $\P(\widehat W^{ \Gamma^N}=  w_i^N)$ are known. The Voronoi quantization  defined by $
x^{(N)}$  induces a lower quantization error but we have no  access to its weights for
numerics. The good news is that $\widetilde X^{ x^{(N)}}$is already rate optimal.
\begin{theorem} (\cite{LUPA3} (2006)) Assume that $b$ is differentiable, $\vartheta$
is positive twice differentiable and that $b'-b\frac{\vartheta'}{\vartheta} -\frac 12
\vartheta\vartheta"$ is bounded. Then
\[
e_{_N}(X,L^2_{_T})\le \|\ X-\widetilde X^{ x^{(N)}}\|_{_2}=O((\log N)^{-\frac 12}).
\]
If furthermore, $\vartheta \ge \varepsilon_0>0$, then $\displaystyle
e_{_N}(X,L^2_{_T})\approx  (\log N)^{-{\frac 12}}$.
\end{theorem}

\noindent {\bf Remarks.} $\bullet$ For some results in the non homogenous case, we refer
to~\cite{LUPA3}. Furthermore, the above estimates still hold true for the $(L^r(\P),
L^p_{_T})$-quantization, $1< r,p<+\infty$ provided $\| |W-\widehat W^{\Gamma^N}|_{L^p_{_T}}\|_r=
O((\log N)^{-\frac 12})$.

\smallskip
\noindent   $\bullet$ This result is closely connected to the Doss-Sussman
approach (see $e.g.$~\cite{DOSS}) and in fact the results can be extended to some
classes multi-dimensional diffusions (whose diffusion coefficient is the inverse of the
gradient of a diffeomorphism) which include several standard multi-dimensional financial
models (including the Black-Scholes model).

\smallskip
\noindent   $\bullet$ A sharp quantization rate $e_{_{N,r}}(X,L^p_{_T})\sim c(\log N)^{-\frac
12}$ for scalar  elliptic diffusions  is established in~\cite{Dereich, Dereichb}   using a non
constructive approach, $1\le p\le\infty$.

\medskip
 \noindent {\sc Example:} Rate optimal product quantization   of the
Ornstein-Uhlenbeck process.
\[
dX_t = -kX_t dt + \vartheta dW_t,\qquad
X_0=x_0.
\]
One solves the non-coupled integral (linear) system
\[
x_{i}(t) = x_0 -k\int_0^t
x_{i}(s)\,ds+  \vartheta
w^{N}_{i}(t),
\]
where $\Gamma^{N}:=\{w^N_1,\ldots,w^N_N\},N\ge 1$ is a {\em rate optimal} sequence of
quantizers
\[
w^{N}_{i}(t)= \sqrt{\frac 2T} \sum_{\ell\ge 1}
\varpi_{i,\ell} \frac{T}{\pi(\ell-1/2)}\sin\left(\pi(\ell-1/2)\frac{t}{T}\right),\quad i\!\in I_N.
\]
If $\Gamma^{N}$ is optimal for $W$ then $\varpi_{i,\ell} := (\beta^N_i)^\ell$, $i=1,\ldots,N$, $1\le
\ell\le d(N)$ with the notations introduced in~(\ref{optizedquantW}). If $\Gamma^{N}$ is an  optimal product quantizer
(and $N_1,\ldots,N_\ell,\ldots$ denote the optimal size allocation), then $\varpi_{i,\ell}
=\alpha^{(N_\ell)}_{i_\ell}$, where
$i:=(i_1,\ldots,i_{\ell},\ldots)\!\in \prod_{\ell\ge 1}
\{1,\ldots,N_\ell\}$. Elementary computations show that
\begin{eqnarray*}
x^N_{i}(t) &  = &  e^{-kt}x_0 + \vartheta
\sum_{\ell\ge 1} \chi^{(N_\ell)}_{i_\ell}
\, \widetilde c_\ell\,  \varphi_\ell(t)\\
\mbox{with } \qquad
\widetilde c_\ell &=& \frac{T^2}{(\pi(\ell-1/2))^2 +(kT)^2}\\
\mbox{and}\hskip 0,4 cm  \varphi_\ell(t)&\!:=\!& \sqrt{\frac
2T}\!\left(\!\frac{\pi}{T}(\ell\!\!-\!\!1/2)
\sin\left(\pi(\ell\!\!-\!\!1/2)\frac{t}{T}\right)
\!+\!k\left(\!\cos\left(\!\pi(\ell\!\!-\!\!1/2)\frac{t}{T}\right)\!\!-\!e^{-kt}\right)\!\right).
\hskip 1 cm
\end{eqnarray*}

\subsection{Multi-dimensional diffusions for Stratanovich SDE's}

The correcting term $-\frac 12 \vartheta\vartheta'$ coming up in the integral equations
 suggest to consider directly some diffusion in the Stratanovich sense
\[
dX_t = b(t,X_t)\,dt +\vartheta(t,X_t)\, \circ\, dW_t \qquad X_0=x_0\!\in\R^{d},\qquad t\!\in [0,T].
\]
(see $e.g.$~\cite{REYO} for an
introduction) where $W= (W^1,\ldots,W^d)$ is a $d$-dimensional standard Brownian Motion.
 
In that framework, we need to introduce the notion of $p$-variation: a continuous function $x:[0,T]\to \R^d$ has
finite
$p$-variations if
\[
 Var_{p, [0,T]}(x):=
\sup\left\{\hskip -0,15 cm\left(\sum_{i=0}^{k-1}|x(t_i)-x(t_{i+1})|^p\hskip -0,13 cm\right)^{\frac
1p}\hskip -0,15 cm , 0\le t_0\le t_1\le \dots \le t_k \le T,\; k\ge 1\hskip -0,1 cm\right\}<+\infty.
\]
Then $d_p(x,x') = |x(0)-x'(0)|+Var_{p, [0,T]}(x-x')$ defines a distance on the set of  functions with finite
$p$-variations. It is classical background that $Var_{p, [0,T]}(W(\omega)) <+\infty$ $\P(d\omega)$-$a.s.$ for every $p>2$.

One way to quantize $W$ at level (at most) $N$ is to quantize each component $W^i$ at level $\lfloor
\sqrt[d]{N}\rfloor$. One shows (see~\cite{LUPA2}) that $\|W - (\widehat W^{1,\lfloor
\sqrt[d]{N}\rfloor},\ldots,\widehat W^{d,\lfloor
\sqrt[d]{N}\rfloor})\|_{_2}=O((\log
N)^{-\frac 12})$.

\bigskip
Let ${\cal C}_b^{r}([0,T]\times\R^d)$ $r>0$, denote the set of $\lfloor r\rfloor$-times differentiable bounded functions
$f:[0,T]\times\R^d\to \R^d $ with bounded partial derivatives up to order $\lfloor r\rfloor$ and whose partial derivatives of
order $\lfloor r\rfloor$ are $(r-\lfloor r\rfloor)$-H\"older.
 
\begin{theorem} (see~\cite{PASE}) Let
$b,\vartheta \!\in {\cal C}_b^{2+\alpha}([0,T]\times\R^d)$ 
$(\alpha>0)$ and let
$\Gamma^N=\{w^N_1,\ldots,w^N_{_N}\}$, $N\ge 1$, be a sequence of   $N$-quantizers of the standard $d$-dimensional Brownian
motion $W$ such that  
$\|W-\widehat W^{\Gamma^N}\|_{_2}\to 0$ as $N\to \infty$. Let
\[
\widetilde X^{x^{(N)}}_t :=\sum_{i=1}^N  x_i^{(N)}\!(t)\mbox{\bf 1}_{ \{\widehat
W= w_i^N\}}
\]
where, for every $i\!\in \{1,\ldots,N\}$, $ x_i^{(N)}$ is solution to
$$
ODE_i\quad \equiv\quad \displaystyle dx_i^{(N)}\!(t)=
b(t,x_i^{(N)}\!(t))dt+\vartheta(t,x_i^{(N)}\!(t))d w_i^N(t),\quad x^{(N)}_i(0)=x.
$$
Then, for every $p\!\in (2,\infty)$,
\[
 Var_{p, [0,T]}(\widetilde X^{x^{(N)}} - X) \stackrel{\P}{\longrightarrow}0\quad \mbox{ as }\quad  N\to \infty.
\]
\end{theorem}

\noindent {\bf Remarks.} $\bullet$  The keys of this results are
the  Kolmogorov criterion,  stationarity (in a slightly extended sense) and the
connection with  rough paths theory (see~\cite{LEJ} for an introduction to rough paths
theory, convergence in $p$-variation, etc).

\noindent $\bullet$  In that general setting we have no convergence  rate although we conjecture that
$\widetilde X^{x^{(N)}}$ remains rate optimal if $\widehat W^{\Gamma^N}$ is.

\noindent $\bullet$ There are also some results about the convergence of stochastic integrals of the form
$\displaystyle \int_0^tg(\widehat W^N_s)\,d \widehat B^N_s\to \int_0^tg(W_s)\circ dB_s$, with some rates of
convergence when $W=B$ or $W$ and $B$ independent (depending on the  regularity of the function $g$, see~\cite{PASE}).
\section{Applications to path-dependent option pricing}\label{pathdeppric}
The  typical functionals $F$ defined on $(L^2_{_{T}},|\,.\,|_{L^2_{_{T}}})$ for which $\E\,(F(W))$ can be
approximated by the  cubature formulae~(\ref{Cub1}),~(\ref{Cub2})  are of the form $\displaystyle
F(\omega):=\varphi\left(\int_0^T f(t,\omega(t)) dt\right)\mbox{\bf 1}_{\{\omega\in {\cal C}([0,T],\R)\}} $  where
$f:[0,T]\times
\R\to
\R$  is {\em locally Lipschitz continuous} in the second variable, namely
\[
\forall\, t\!\in [0,T],\; \forall\, u,v\!\in \R,\; |f(t,u)-f(t,v)|\le C_f |u-v|(1+g(|u|)+g(|v|))
\]
(with $g:\R_+\to\R_+$ is increasing, convex and $g(\sup_{t\in[0,T]}|W_t|)\!\in L^2(\P)$) and $\varphi:\R\to \R$ is
Lipschitz continuous. One could consider for $\omega$ some c\`adl\`ag functions as well. A classical example is the
Asian payoff in a  Black-Scholes  model
\[
F(\omega) =\exp(-rT)\left( \frac 1T \int_0^T
s_0\exp(\sigma\omega(t)+(r-\sigma^2/2)t)dt-K\right)_+.
\]
\subsection{Numerical integration (II): $\log$-Romberg
extrapolation}

Let  $F: L^2_{_T} \longrightarrow \R$ be a $3$ times
$|\,.\,|_{L^2_{_T}}$-differentiable functional with bounded differentials. Assume  $\widehat W ^{(N)}$, $N\ge
1$,  is a sequence of a  rate-optimal stationary quantizations of the standard Brownian motion $W$. Assume
furthermore that
\begin{equation}\label{Conj1}
\E\left(D^2F(\widehat W ^{(N)}).(W-\widehat W^{(N)})^{\otimes
2}\right)\sim \frac{c}{\log N}\quad \mbox{ as }\quad N\to
\infty
\end{equation}
and
\begin{equation}\label{Conj2}
  \E\,|W-\widehat W^{(N)}|_{L^2_{_T}}^3 =
O\left((\log N)^{-\frac 32}\right).
\end{equation}
Then,  a higher order Taylor expansion  yields
\begin{eqnarray*}
 \quad F(W)&=& F(\widehat W ^{(N)}) + DF(\widehat W ^{(N)}).(W-\widehat W ^{(N)})
 +\frac 12 D^2F(\widehat W ^{(N)}).(W-\widehat W^{(N)})^{\otimes 2} \\
&& + \frac 16 D^2(\zeta).(W-\widehat W ^{(N)})^{\otimes
3}, \qquad \zeta\!\in (\widehat W ^{(N)},W),\\
\E \,F(W)&=&\E F(\widehat W ^{(N)})  + \frac{c}{2\log N} +o\left((\log N)^{-\frac
32+\varepsilon}\right).\hskip 2 cm
\end{eqnarray*}
Then, one can design a $\log$-Romberg extrapolation by considering $N,\,N'$, $N < N'$ ($e.g.$ $N'\approx
4\,N$), so that
\[
\hskip -0.25 cm     \E(F(W)) =\frac{ \log
N'\!\times\!\E(F(\widehat{W}^{(N')}))-\log
N'\!\times\!\E(F(\widehat{W}^{(N)})) }{\log N'-\log N} + o\left((\log
N)^{-\frac 32+\varepsilon}\right).
\]
For practical implementation, it is suggested in~\cite{WIL} to replace $\log N$ by the  more
consistent  ``estimator"
$\|W-\widehat W^{(N)}\|_{_2}^{-2}$.

\smallskip In fact Assumption~(\ref{Conj1}) holds true for optimal product quantization when $F$ is
polynomial function $F$, $d^0F=2$. Assumption~(\ref{Conj2}) holds true in that case as well
(see~\cite{GRLUPA3}). As concerns  optimal quantization, these statements are still conjectures. However, given
that
$\widehat W$ and
$W-\widehat W$ are  independent (see~\cite{LUPA1}),~(\ref{Conj1}) is equivalent to the simple
case where
$D^2F(\widehat W ^{(N)})$ is constant.

Note that the above extrapolation or some variants  can be implemented with other stochastic processes in accordance
with the rate of convergence of the quantization error.

\subsection{Asian option pricing in a Heston stochastic volatility model}
In this section, we will price an Asian call option in a Heston stochastic volatility model using some optimal
(at least optimized) functional quantization of the two Brownian motions that drive the diffusion. This model
has already been considered   in~\cite{PAPR2} in which  functional
quantization was implemented for  the first time with some product quantizations of the Brownian motions.
 The Heston stochastic volatility model was introduced
in~\cite{HES} to model stock price dynamics. Its popularity partly comes from the existence of
semi-closed forms for vanilla European options, based on inverse Fourier transform and from its ability to
reproduce some skewness shape of the implied volatility surface. We consider it under its risk-neutral
probability measure.
\begin{eqnarray*}
dS_t&=& S_t(r\,dt +\sqrt{v_t}dW^1_t),\qquad S_0=s_0>0, \quad \mbox{(risky asset)} \\
\nonumber dv_t&=&k(a-v_t)dt +\vartheta \sqrt{v_t}\,dW^2_t,\; v_0> 0 \; \mbox{with }\; d\!<\!\!W^1,W^2\!\!>_t
=\rho \,dt,\;\rho\!\in[-1,1].
\end{eqnarray*}
where $\vartheta, k, a$ such that
$\vartheta^2/(4ak)< 1$. We  consider the Asian Call payoff with maturity $T$ and strike $K$. No  closed form  is available for its premium
\[
{\rm AsCall}^{Hest}= e^{-rT}  \E\left(\frac 1T \int_0^T S_s ds
-K\right)^+.
\]
We briefly recall how to proceed (see~\cite{PAPR2} for details): first, one projects
$W^1$ on $W^2$ so that $W^1 =\rho W^2+\sqrt{1-\rho^2}\,\widetilde W^1$ and
\begin{eqnarray*}
S_t&= &s_0\exp{\left(\hskip -0.1cm(r-\frac{1}{2}\bar
v_t)t+\rho\int_0^t\sqrt{v_s}dW^2_s\right)}\exp{\left(\sqrt{1-\rho^2}\int_0^t\sqrt{v_s}d\widetilde{W}^1_s\right)}
\\
&=& s_0\exp{\left(\hskip -0.1cmt\left(\hskip -0.1cm(r-\frac{\rho a k}{\vartheta})+\bar
v_t(\frac{\rho k}{\vartheta}-\frac{1}{2})\hskip -0.1cm\right)+
\frac{\rho}{\vartheta}(
v_t-v_0)\hskip -0.1cm\right)}\exp{\left(\hskip -0.1cm\sqrt{1-\rho^2}\int_0^t\sqrt{v_s}d\widetilde{W}^1_s\hskip
-0.1cm\right)}.
\end{eqnarray*}
The chaining rule for conditional expectations yields
\[
{\rm AsCall}^{Hest}(s_0,K)= e^{-rT} \E\hskip -0.1  cm\left(\hskip -0.1  cm \E\hskip -0.1cm\left(\hskip -0.1
cm\left(\frac 1T
\int_0^T S_sds -K\hskip -0.1  cm\right)^+\hskip -0.15  cm|\sigma(W^2_t, 0\le t\le T)\hskip -0.1  cm\right)\hskip
-0.1  cm\right).
\]
Combining these two expressions and using  that $\widetilde W^1$ and $W^2$ are independent show that
${\rm AsCall}^{Hest}(s_0,K)$ is a functional of $(\widetilde W^1_t, v_t)$ (as concerns the squared volatility
process $v$, only $v_{_T}$ and $\int_0^T v_sds$ are involved).

\smallskip
Let $\Gamma^N=\{w^N_1,\ldots,w^N_{_N}\}$ be an $N$-quantizer of the Brownian motion. One solves
for $ i=1,\ldots,N$, the differential equations  for $(v_t)$
\begin{equation}\label{ODEi}
dy_i(t) =  k\left(a-y_i(t)-\frac{\vartheta^2}{4k}\right) dt +
\vartheta \sqrt{y_i(t)}\,dw^{N}_i(t), \;y_i(0)=v_0,
\end{equation}
 using
$e.g.$ a Runge-Kuta scheme. Let $y^{n,N}_i$ denote the approximation of $y_i$
resulting from the resolution of the above $ODE_i$ ($1/n$ is the time discretization parameter of the
scheme).
Set the  (non-Voronoi) $N$-quantization of $(v_t,S_t)$ by
\begin{eqnarray}\label{QuantifieurV1}
\widetilde v_t^{n,N}&=& \sum_{i} y^{n,N}_{i}(t)\mbox{\bf 1}_{C_{i}(\Gamma^N)}(W^2)\\
\label{chaineAsi2}
  \widetilde S^{n,N}_t &=&\sum_{ 1\le i, j\le N}
 s^{n,N}_{ i, j}(t) \mbox{\bf 1}_{C_i(\Gamma^N)}(\widetilde W^1)\mbox{\bf
1}_{C_j(\Gamma^N)}( W^2)  \\
 \nonumber \mbox{with  }\;s^{n,N}_{ i, j}(t) &=&s_0\exp{\left(\!t\!\left(\!(r-\frac{\rho a
k}{\vartheta})+\overline y^{n,N}_{ j}(t)(\frac{\rho
k}{\vartheta}-\frac{1}{2})\!\right)+ \frac{\rho}{\vartheta}(
y^{n,N}_{j}(t)-v_0)\!\right)}\quad\\
 \nonumber  &&\times \exp{\left(\sqrt{1-\rho^2} \int_0^t\sqrt{y^{n,N}_{j}(s)}\,dw^N_{
i}(s)\right)}\\
 \nonumber \mbox{and }\quad \overline y^{n,N}_{ j}(t)&=&\int_0^ty^{n,N}_{ j}(s)\,ds.
\end{eqnarray}
Note this formula requires the computation of a quantized stochastic integral
$\displaystyle \int_0^t\sqrt{y^{n,N}_{j}}(s)dw^N_{ i}(s)$ (which corresponds to the independent
case).

\noindent The weights of the product cells $\{\widetilde W^1\!\in C_i(\Gamma^N),\, W^2\!\in C_j(\Gamma^N)\}$ is
given by
\[
\P(\widetilde W^1\!\in C_i(w^N),\, W^2\!\in C_j(w^N))= \P(\widetilde W^1\!\in C_i(\Gamma^N))\P(W^2\!\in
C_j(\Gamma^N))
\]
owing to the independence. For practical implementations different
sizes of quantizers can be considered to quantize $\widetilde W^1$
and $W^2$.

\medskip We follow the guidelines of the methodology introduced in~\cite{PAPR2}: we     compute the {\em crude}
quantized premium for two sizes $N$ and $N'$, then  proceed   a space Romberg $\log$-extrapolation. Finally, we
make a $K$-linear interpolation   based on the (Asian) forward moneyness
$s_0e^{rT}\frac{1-e^{-rT}}{rT}\approx s_0e^{rT}$ (like in~\cite{PAPR2})
and the  Asian Call-Put parity  formula
\[
{\rm AsianCall}^{Hest}(s_0,K)={\rm AsianPut}^{Hest}(s_0,K)+
s_0\frac{1-e^{-rT}}{rT}-Ke^{-rT}.
\]
The  {\em anchor strikes} $K_{\min}$ and $K_{\max}$ of the
extrapolation are chosen symmetric with respect to the forward
moneyness. At $K_{\max}$, the Call is deep out-of-the-money: one
uses the  Romberg extrapolated $FQ$ computation; at $K_{\min}$ the
Call is deep in-the-money: on computes the Call by parity.  In
between, one proceeds a linear interpolation in $K$ (which yields
the best results, compared to other extrapolations like the
quadratic regression approach).
\begin{figure}
\centering
\begin{tabular}{cc}
\includegraphics[width=6cm,height = 5cm]{400_6.eps}&\includegraphics[width=6cm,height = 4.5cm]{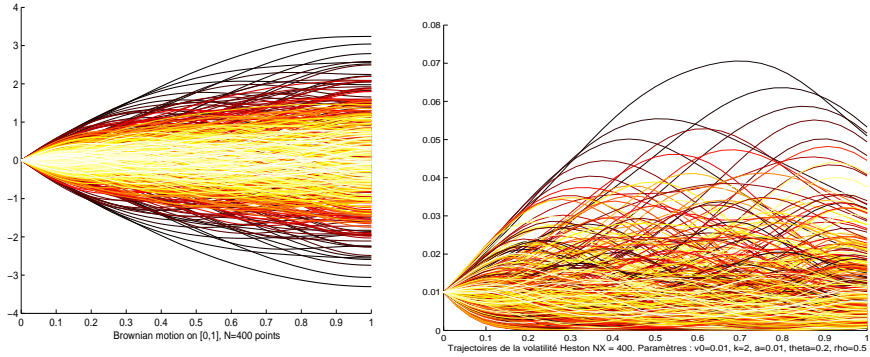}
\end{tabular}
\caption{\em $N$-quantizer  of the Heston squared volatility process $(v_t)$ ($N=400$) resulting from an (optimized)
$N$-quantizer of $W$.}
\end{figure}

\smallskip
$\circ$ {\em Parameters of the Heston model}: $s_0=100$, $k=2$, $a=
0.01$, $\rho = 0.5$, $v_0 =10\%$, $\vartheta=20\%$.

\smallskip
$\circ$ {\em Parameters of the option portfolio}: $T=1$, $K=99,\cdots,
111$ (13 strikes).

\smallskip
$\circ$ {\em The reference price} has been   computed by a $10^8$
trial Monte Carlo simulation  (including a time Romberg
extrapolation of the Euler scheme with $2n=256$).
 
\smallskip
$\circ$ {\em The differential equations~(\ref{ODEi})} are solved with the
parameters of the quantization cubature formulae $\Delta t = 1/32$,
with couples of quantization levels $(N,M)= (400,100)$,
$(1000,100)$, $(3200,400)$.

\begin{figure}
\centering
\hskip -0.25 cm
\includegraphics[width=12cm,height= 5cm]{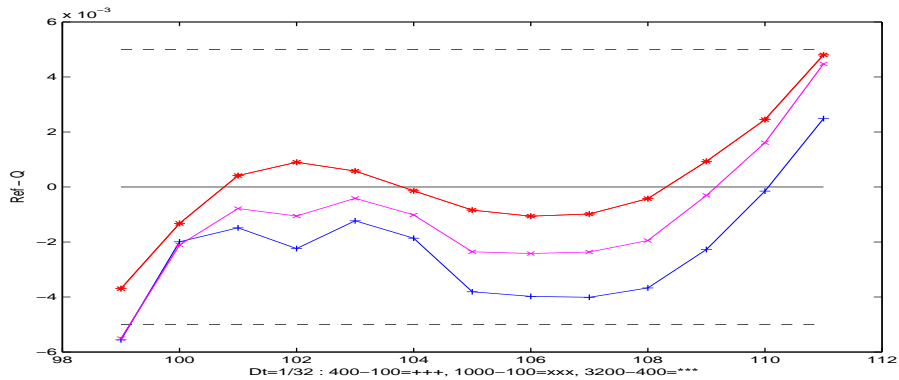}
\caption{\em Quantized diffusions based on optimal functional
quantization: Pricing by $K$-Interpolated-$\log$-Romberg extrapolated-$FQ$
prices as a function of $K$: absolute error with
$(N,M)\!=\!(400,100)$, $(N,M)\!=\!(1000,100)$,
$(N,M)\!=\!(3200,400)$. $T\!=\!1$, $s_0\!=\!50$, $K\!\in\{99,\ldots,
111\}$.   $k\!=\!2$, $a\!=\! 0.01$, $\rho \!=\! 0.5$,
$\vartheta\!=\!0.1$.}
\end{figure}

\begin{figure}
\centering
\hskip -0.25 cm
\includegraphics[width= 5cm,height= 12cm, angle =270]{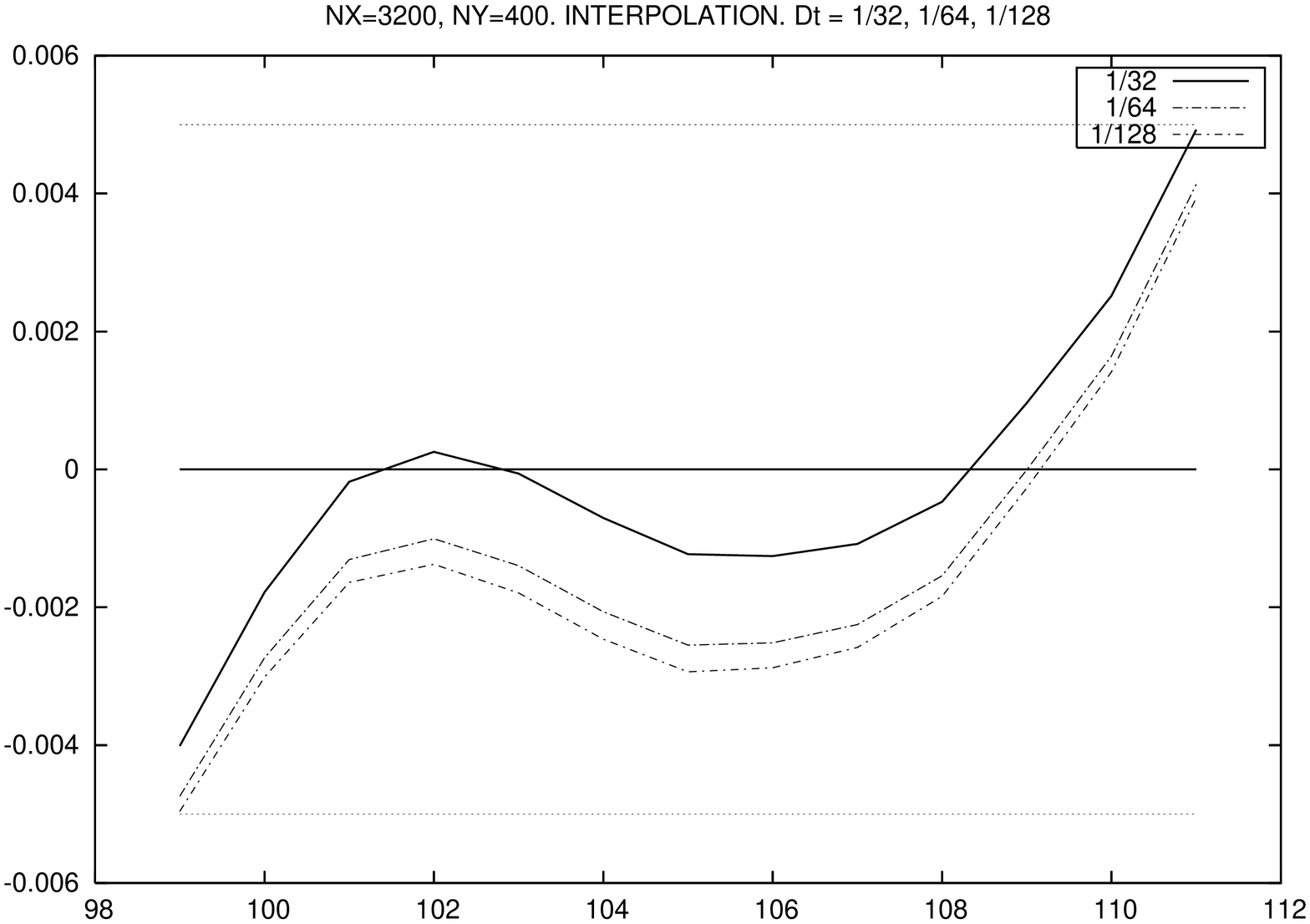}
\caption{\em Quantized diffusions based on optimal functional
quantization: Pricing by $K$-Interpolated-$\log$-Romberg extrapolated-$FQ$
price as a function of $K$: convergence as $\Delta t\to 0$ with
$(N,M)=(3200,400)$ (absolute error). $T\!=\!1$, $s_0\!=\!50$,
$K\!\in\{99,\ldots, 111\}$.   $k\!=\!2$, $a\!=\! 0.01$, $\rho \!=\!
0.5$, $\vartheta\!=\!0.1$.}
\end{figure}

Functional Quantization can compute  a
whole vector (more than $10$) option premia for the Asian option in
the Heston model  with {\bf  $1$ cent accuracy
in less than $1$ second} (implementation  in $C$ on a $2.5$ $GHz$ processor).

\smallskip
Further numerical tests carried out or in progress  with the $B$-$S$
model and with the  $SABR$ model (Asian, vanilla European options)
show the same efficiency. Furthermore, recent attempt to quantize
the volatility process and the asset dynamics at different level of
quantizations seem very promising in two directions: reduction of
the computation time  and increase of the robustness of the method
to parameter change.

\subsection{Comparison: optimized quantization $vs$  (optimal)
product quantization}
\begin{figure}\label{fig:AsCall_B}
\centering
\begin{tabular}{c}
$$
\hskip -0.25 cm \includegraphics[width=12cm,height
=5cm]{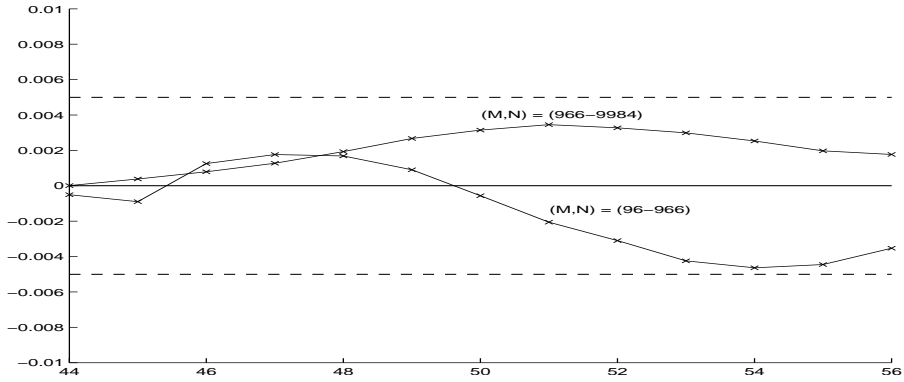}
$$
\end{tabular}
\caption{\em  Quantized diffusions based on optimal product
quantization:
Pricing by $K$-linear
interpolation of Romberg $\log$-extrapolations as un function of $K$ (absolute error)
with $(M,N)$= $(96,966)$,
$(966,9984)$.  $T\!=\!1$, $s_0\!=\!50$,   $k\!=\!2$, $a\!=\! 0.01$, $\rho \!=\!0.5$, $\vartheta\!=\!0.1$.
$K\!\in\{44,\ldots, 56\}$.}
\end{figure}
The comparison is balanced and  probably  needs  some further {\em in situ} experiments since it may depend on the
modes of the computation. However, it seems that product quantizers (as those implemented in~\cite{PAPR2}) are 
from $2$   up to $4$ times less efficient  than optimal quantizers within our range of application (small values of
$N$). On the other hand, the design  of product quantizer  from $1$-dim scalar quantizers is easy and can be made
from some light elementary ``bricks" (the scalar quantizer up to $N=35$ and the  optimal allocation rules).  
Thus, the whole set of data needed to design all  optimal product quantizers up to $N=10\,000$ is approximately 
$500$ $KB$ whereas one optimal quantizer with size $10\, 000  \approx 1~MB$\dots 
\section{Universal quantization rate and mean regularity}\label{FQReg}
The following theorem points out the connection between functional
quantization rate and mean regularity of $t\mapsto X_t$ from $[0,T]$
to $L^r(\P)$.
\begin{theorem} (\cite{LUPA4} (2005)) Let
$X=(X_t)_{t\in[0,T]}$ be a stochastic process. If there is $r^*\!\in (0,\infty)$ and $a\!\in (0,1]$ such that
\[
X_0 \in L^{r^*}(\P),\quad \|X_t-X_s\|_{L^{r^*}(\P)}\le C_X|t-s|^{a},
\]
for some positive real constant $C_X>0$, then
$$
\displaystyle \forall\,p,r\!\in (0,r^*),\quad
e_{N,r}(X,L^p_{_T}) = O((\log N)^{-a}).
$$
\end{theorem}
The  proof is based on a  constructive approach which involves the Haar basis (instead of $K$-$L$
basis), the  non asymptotic version Zador Theorem and product functional quantization. Roughly speaking, we use the
unconditionality of the Haar basis in every $L^p_{_T}$  (when $1\!< \!p<\infty$) and its
wavelet feature $i.e.$ its ability to ``code" the path regularity of a function on the decay
rate of its coordinates.

\medskip
\noindent {\sc  Examples (see~\cite{LUPA4}):}   $\bullet$
$d$-dimensional It\^o processes (includes  $d$-dim
 diffusions with sublinear coefficients) with $a= 1/2$.

\noindent $\bullet$  General  L\'evy process $X$ with L\'evy measure $\nu$ with
square integrable big jumps. If $X$ has a Brownian component, then $a=2$, otherwise if
 $\beta(X)>0$ where $\beta(X) :=\inf\left\{\theta\,:\, \int|y|^\theta\nu(dy)
\!<\!+\infty\right\}\!\in (0,2)$ (Blumenthal-Getoor index of $X$), then $a=   \beta^*(X)$. This rate is the {\em
exact rate} $i.e.$
\[
e_{N,r}(X,L^p_{_T}) \approx (\log N)^{-a}
\]
for many classes of L\'evy processes like symmetric stable processes, L\'evy processes having a Brownian component, etc (see~\cite{LUPA4} for further examples).

\noindent$\bullet$ When $X$ is a compound Poisson processes, then  $\beta(X)=0$ and one shows, still
with   constructive methods,   that
 \[
 e_N(X) = O(e^{-(\log N)^{\vartheta}}),\qquad\vartheta \!\in (0,1),
 \]
 which is in-between the finite and infinite dimensional settings.

\section{About lower bounds}\label{LowBounds} In this overview, we gave no clue toward lower bounds although most of the rates
we mentioned are either exact ($\approx$) or sharp ($\sim$) (we tried to emphasize the numerical
aspects). Several approaches can be developed to get some lower bounds. Historically,
the first one was to rely on subadditivity property of the
quantization error derived from self-similarity of the distribution:
this works with the uniform distribution over $[0,1]^d$ but also in
an infinite dimensional framework (see~$e.g.$~\cite{Dereichetal} for
the fractional Brownian motion).

A second approach consists in pointing out the  connection with the Shannon-Kolmogorov entropy
(see~$e.g.$~\cite{LUPA1}) using that the entropy of a random variable taking at most $N$ values is at most $\log
N$.

A third connection can be made with small deviation theory (see~\cite{DEFEMASC},~\cite{GRLUPA1} and~\cite{LUPA4}). Thus,
 in \cite{GRLUPA1}, a connection is established between (functional) quantization and small ball deviation for Gaussian processes. In particular
this approach provides a method to derive a lower bound for the
quantization rate from some upper
bound for the small deviation problem. A careful reading of   the proof of Theorem~1.2 in~\cite{GRLUPA1} shows
that this small deviation lower bound holds for any {\em unimodal} (w.r.t. $0$) non zero process. To be precise:
 assume that  $\P_{_X}$ is $L^p_{_T}$-unimodal $i.e.$ there exists a real
$\varepsilon_0 > 0$ such that
\[
\forall\, x\!\in L^p_{_T},\;\forall\, \varepsilon \!\in (0,\varepsilon_0],\qquad  \P(|X-x|_{L^p_{_T}}\le \varepsilon)\le
\P(|X|_{L^p_{_T}}\le \varepsilon).
\]
For centered Gaussian processes (or processes ``subordinated" to
Gaussian processes) this follows from the Anderson  Inequality (when $p\ge 1$). If
\[
G(-\log(\P(|X|_{L^p_{_T}}\le \varepsilon)))= \Omega(1/\varepsilon) \quad \mbox{as} \quad \varepsilon\to 0
\]
for some increasing unbounded function $G:(0,\infty)\to (0,\infty)$, then
\begin{equation}\label{sblowerbound}
\forall\, c>1,\quad \liminf_N G(\log(cN)) e_{_{N,r}}(X, L^p_{_T})>0,\qquad r\!\in(0,\infty).
\end{equation}
This approach is   efficient in the non quadratic case as emphasized  in~\cite{LUPA4} where several universal bounds are shown
to be optimal using this approach.

\medskip
\noindent {\sc Acknowledgement.} I  thank S.
Graf, H. Luschgy J. Printems and B. Wilbertz for all the  fruitful
discussions and collaborations we have about functional quantization.


\begin{thebibliography}{0}
\bibitem{ABA1}{\sc Abaya, E.F. and Wise, G.L.} (1982). On the existence of optimal
quantizers. {\em IEEE Trans. Inform. Theory}, {\bf  28}, 937-940

\bibitem{ABA2}
{\sc Abaya, E.F. and Wise, G.L.} (1984). Some remarks on the
existence of optimal quantizers. {\em Statistics and Probab.
Letters}, {\bf 2}, 349-351.


\bibitem{BMP}{\sc Benveniste, A., M\'etivier, M. and Priouret, P.} (1990).
{\em Adaptive algorithms and stochastic approximations},
Translated from the French by Stephen S. Wilson. Applications of Mathematics {\bf 22}, Springer-Verlag, Berlin, 365 pp.

\bibitem{BOLE} {\sc N. Bouleau, D. L\'epingle} (1994).  {\em
Numerical methods for stochastic processes}, Wiley Series in Probability and Mathematical Statistics: Applied Probability and
Statistics. A Wiley-Interscience Publication. John Wiley
\& Sons, Inc.,  New York, 359 pp. ISBN: 0-471-54641-0.


\bibitem{BUWI1} {\sc Bucklew, J.A. and Wise, G.L.} (1982). Multidimensional
asymptotic quantization theory with $r^{th}$ power distortion. {\em
IEEE Trans. Inform. Theory}, {\bf  28}(2), 239-247.

\bibitem{COH} {\sc Cohort, P.} (1998).  A geometric method for  uniqueness of
locally optimal quantizer. Pre-print LPMA-464 and Ph.D.
Thesis, {\em Sur quelques probl\`emes de quantification}, 2000,
Univ. Paris 6.

\bibitem{CUMA} {\sc Cuesta-Albertos, J.A., Matr{\'a}n, C.} (1988). The strong law of large numbers for $k$-means and best
possible nets of Banach valued random variables, {\em Probab. Theory Rel. Fields }{\bf 78},  523-534.

\bibitem{DELA} {\sc Delattre, S., Fort, J.-C. and Pag\`es, G.} (2004). Local
distortion and $\mu$-mass of the cells of one dimensional
asymptotically optimal quantizers, {\em Communications in
Statistics}, {\bf 33}(5),  1087-1118.

\bibitem{DEFEMASC} {\sc Dereich, S., Fehringer, F., Matoussi, A. and Scheutzow, M.}
(2003). On the link between small ball probabilities and the
quantization problem for Gaussian measures on Banach spaces, {\em J.
Theoretical Probab.}, {\bf 16},   pp.249-265.
%
\bibitem{Dereich} {\sc Dereich, S.} (2005). The coding complexity of diffusion processes under $L^p[0,1]$-norm
distortion, pre-print.
%
\bibitem{Dereichb} {\sc Dereich, S.} (2005).  The coding complexity of diffusion processes under supremum norm distortion,
pre-print.
%
\bibitem{Dereichetal}
{\sc Dereich, S., Scheutzow, M.} (2006). High resolution quantization and
entropy coding for fractional Brownian motions, {\em Electron. J. Probab.},  {\bf 11},  700-722.

\bibitem{DOSS} {\sc Doss  H.} (1977). Liens entre \'equations 
diff\'erentielles stochastiques et ordinaires, {\em Ann. I.H.P.}, 
section B, {\bf 13}(2),  99-125.

\bibitem{FLE}
{\sc Fleischer, P.E.} (1964). Sufficient conditions for achieving
minimum distortion in a quantizer. {\em IEEE Int. Conv. Rec.}, part
I, 104-111.

\bibitem{GEGR} {\sc Gersho, A. and Gray, R.M.} (1992). {\em
Vector Quantization and Signal Compression}. Kluwer, Boston.

\bibitem{FOPA}
{\sc  Fort, J.-C. and Pag\`es, G.} (2004). Asymptotics
of optimal quantizers for some scalar distributions,  {\em Journal of  Computational
and Applied  Mathematics}, {\bf 146}, 253-275, 2002.

\bibitem{GRLU1}
{\sc Graf, S. and Luschgy, H.} (2000). {\em Foundations of
Quantization for Probability Distributions}. Lect. Notes in Math.
1730, Springer, Berlin, 230p.

\bibitem{GRLU2}
{\sc Graf, S. and Luschgy, H.} (2005). The point density measure in the quantization of self-similar  probabilities. {\em
Math. Proc. Cambridge Phil. Soc.}. {\bf 138}, 513-531.

\bibitem{GRLUPA1} {\sc Graf, S., Luschgy H. and Pag\`es, G.} (2003). Functional
quantization and small ball probabilities for Gaussian processes,  {\em J. Theoret. Probab.}, {\bf 16}(4), 1047-1062.

\bibitem{GRLUPA2}
{\sc Graf,  S., Luschgy,  H., Pag\`es,  G.} (2007). Optimal quantizers for Radon random vectors in a Banach space,
 {\em J. of Approximation Theory}, {\bf 144}, 27-53.

\bibitem{GRLUPA3}  {\sc Graf, S., Luschgy,  H. and  Pag\`es,  G.} (2006). Distortion mismatch in the
quantization of probability measures,
to appear in {\em ESAIM P\&S}.

\bibitem{HES} {\sc Heston, S.L.} (1993). A closed-form solution for
options with stochastic volatility with applications to bond and
currency options, {\em
The review of Financial Studies}, {\bf 6}(2), 327-343.

\bibitem{KIE} {\sc Kieffer, J.C.} (1983). Uniqueness of locally optimal quantizer for
$\log$-concave density and convex error weighting functions,  {\em
IEEE Trans. Inform. Theory}, {\bf  29}, 42-47.

\bibitem{KIE2} {\sc Kieffer, J.C.} (1982). Exponential rate of convergence for Lloyd's Method~I,  {\em
IEEE Trans. Inform. Theory}, {\bf  28}(2), 205-210.

\bibitem{KUYI} {\sc Kushner, H. J., Yin, G. G.} (2003). {\em Stochastic approximation and recursive algorithms and
applications}. Second edition. Applications of Mathematics {\bf 35}. Stochastic Modelling and Applied Probability.
Springer-Verlag, New York, 474p.

\bibitem{LAPA} {\sc Lamberton, D. and Pag\`es, G.} (1996). On the critical points of the
$1$-dimensional Competitive Learning Vector Quantization Algorithm.
{\em Proceedings of the ESANN'96}, (ed. M. Verleysen), Editions D
Facto, Bruxelles, 97-106.

\bibitem{LAPASA}{\sc Lapeyre, B., Sab, K. and  Pag\`es, G.} (1990). Sequences with low discrepancy. Generalization and
application to Robbins-Monro algorithm, {\em Statistics}, {\bf 21}(2), 251-272.

\bibitem{LEJ} {\sc Lejay, A.} (2003). An introduction to rough paths, {\em
S\'eminaire de Probabilit\'es XXXVII}, {Lecture Notes in Mathematics 1832}, Stringer, Berlin, 1-59.

\bibitem{LUPA1}  {\sc Luschgy, H., Pag\`es, G.} (2002).  Functional quantization of Gaussian processes,  {\em  Journal of
Functional Analysis},  {\bf 196}(2), 486-531.

\bibitem{LUPA2}  {\sc Luschgy, H., Pag\`es, G.} (2004). Sharp asymptotics of the functional quantization problem
for Gaussian processes,  {\em  The Annals of
Probability},  {\bf 32}(2), 1574-1599.

\bibitem{LUPA3} {\sc Luschgy, H., Pag\`es, G.} (2006). Functional quantization of a class of Brownian
diffusions: A constructive approach,  {\em Stochastic Processes and Applications},  {\bf 116},
310-336.

\bibitem{LUPA3.5} {\sc Luschgy, H., Pag\`es, G.} (2005). High-resolution product quantization for Gaussian processes under
sup-norm distortion, pre-pub LPMA-1029,
forthcoming in {\em Bernoulli}.

\bibitem{LUPA4} {\sc Luschgy, H., Pag\`es, G.} (2006). Functional Quantization Rate and mean
regularity of processes with an application to L\'evy Processes, pre-print
LPMA-1048.

\bibitem{LUPASpec} {\sc Luschgy, H., Pag\`es, G.} (2007). Expansion of Gaussian processes and Hilbert
frames, technical report.

\bibitem{LUPAWI} {\sc Luschgy, H., Pag\`es, G. and Wilbertz, B.} (2007). Asymptotically optimal quantization
schemes for Gaussian processes, in progress.

\bibitem{MRBenH} {\sc Mrad, M., Ben Hamida, S.} (2006). Optimal Quantization: Evolutionary Algorithm $vs$ Stochastic
Gradient,   {\em Proceedings of the 9th Joint Conference on Information Sciences}.

\bibitem{NEW} {\sc Newman, D.J.} (1982). The Hexagon
Theorem.  {\em IEEE Trans. Inform. Theory}, {\bf  28}, 137-138.

\bibitem{NIE} {\sc H. Niederreiter} (1992) {\em Random Number Generation and Quasi-Monte Carlo Methods}, CBMS-NSF regional conference series in Applied mathematics, SIAM, Philadelphia.

\bibitem{PAG0}
{\sc Pag\`es, G.} (1993). Voronoi tessellation, space quantization
algorithm and numerical integration. {\em Proceedings of the
ESANN'93},  M. Verleysen Ed., Editions D Facto, Bruxelles, 221-228.


\bibitem{PAG1} {\sc Pag\`es, G.} (1997). A space vector quantization method for numerical integration,   {\em J.
Computational and Applied Mathematics}, {\bf 89}, 1-38.


\bibitem{PAG0.5} {\sc Pag\`es, G.} (2000). Functional quantization: a first approach, pre-print CMP12-04-00, Univ.
Paris 12.

\bibitem{PAPHPR}
{\sc Pag\`es, G., Pham, H. and Printems, J.} (2003). Optimal
quantization methods and applications to numerical methods in
finance. {\em Handbook of Computational and Numerical Methods in Finance}, S.T. Rachev
ed., Birkh{\"a}user, Boston, 429p.


\bibitem{PAPR1} {\sc  Pag\`es, G., Printems, J.} (2003). Optimal
quadratic quantization for numerics: the Gaussian case, {\em Monte
Carlo Methods and Appl.},  {\bf 9}(2), 135-165.


\bibitem{PAPR2} {\sc Pag\`es, G., Printems, J.} (2005). Functional
quantization for numerics with an application to option pricing,
{\em Monte Carlo Methods and Appl.},  {\bf 11}(4), 407-446.

\bibitem{Website} {\sc Pag\`es, G., Printems, J.} (2005). Website devoted to vector and functional optimal
quantization: {\tt www.quantize.maths-fi.com}.

\bibitem{PASE} {\sc Pag\`es, G., Sellami, A.} (2007). Convergence of multi-dimensional quantized SDE's. In progress.

\bibitem{PAXI}  {\sc G. Pag\`es, Y.J. Xiao} (1988) Sequences with low discrepancy and pseudo-random numbers: theoretical results and numerical tests, {\em  J. of Statist. Comput. Simul.}, {\bf 56}, 163-188.


\bibitem{PAR} {\sc P{\"a}rna,  K.} (1990).  On the existence and weak convergence of
$k$-centers in Banach spaces,  {\em Tartu {\"U}likooli Toimetised}, {\bf 893},  17-287.

\bibitem{POL} {\sc Pollard, D.} (1982). Quantization and the method of $k$-means. {\em IEEE Trans.
Inform. Theory}, {\bf 28}(2), 199-205.

\bibitem{PRO} {\sc P.D. Proinov} (1988). Discrepancy and integration of continuous functions, {\em J. of Approx. Theory}, {\bf 52}, 121-131.

 
\bibitem{REYO}  {\sc Revuz, D., Yor, M.} (1999). {\em Continuous martingales and Brownian motion}, Third edition.
Grundlehren der Mathematischen Wissenschaften [Fundamental Principles of Mathematical Sciences], 293, Springer-Verlag,
Berlin, 1999, 602 p.

\bibitem{ROT}  {\sc K.F. Roth} (1954). On irregularities of distributions, {\em Mathematika}, {\bf 1}, 73-79.


\bibitem{TAKI} {\sc Tarpey, T., Kinateder, K.K.J.} (2003). Clustering functional data, {\em J. Classification}, {\bf 20},
93-114.

\bibitem{TAPEOG} {\sc Tarpey, T., Petkova, E., Ogden, R.T.} (2003). Profiling Placebo responders by self-consistent
partitioning of functional data, {\em J. Amer. Statist. Association}, {\bf 98}, 850-858.

\bibitem{TRU}
{\sc Trushkin, A.V.} (1982). Sufficient conditions for uniqueness
of a locally optimal quantizer for a class of convex error weighting
functions, {\em IEEE Trans. Inform. Theory}, {\bf 28}(2), 187-198.

\bibitem{WIL} {\sc Wilbertz, B.} (2005). Computational aspects of functional quantization for Gaussian measures
and applications, diploma thesis, Univ. Trier.

\bibitem{ZAD0}
{\sc Zador, P.L.} (1963). Development and evaluation of procedures
for quantizing multivariate distributions. Ph.D. dissertation,
Stanford Univ.

\bibitem{ZAD1}
{\sc Zador, P.L.} (1982). Asymptotic quantization error of
continuous signals and the quantization dimension. {\em IEEE Trans.
Inform. Theory}, {\bf 28}(2), 139-149.



\end{thebibliography}
\end{document}